\documentclass[a4paper,12pt]{article}
\usepackage{amsmath,amssymb,graphics}
\newcommand{\field}[1]{\mathbb{#1}}
\usepackage[left=2.4cm,top=3.0cm,bottom=3.1cm,right=2.4cm,head=0cm,foot=0.7cm]{geometry}
\begin{document}
\newtheorem{definition}{Definition}
\newtheorem{lemma}{Lemma}
\newtheorem{proposition}{Proposition}
\newtheorem{theorem}{Theorem}
\newtheorem{remark}{Remark}
\newtheorem{corollary}{Corollary}
\begin{center}
\Large{On the Observational Equivalence of Continuous-Time Deterministic and Indeterministic Descriptions\\}\normalsize
\begin{verbatim}
\end{verbatim}
\large{Charlotte Werndl, c.s.werndl@lse.ac.uk\\
Department of Philosophy, Logic and Scientific Method\\
London School of Economics and Political Science}
\begin{verbatim}
\end{verbatim}
\normalsize{This is a pre-copyedited, author-produced PDF of an article accepted for publication in the European Journal for Philosophy of Science. The definitive publisher-authenticated version ``C. Werndl (2011), On the Observational Equivalence of Continuous-Time Deterministic and Indeterministic Descriptions, European Journal for Philosophy of Science 1 (2), 193-225'' is available online at: http://link.springer.com/article/10.1007\%2Fs13194-010-0011-5.}
\begin{verbatim}
\end{verbatim}
\end{center}
Abstract. This paper presents and philosophically assesses three types of results on the observational equivalence of continuous-time measure-theoretic deterministic and indeterministic descriptions.  The first results establish observational equivalence to abstract mathematical descriptions. The second results are stronger because they show observational equivalence between deterministic and indeterministic descriptions found in science. Here I also discuss Kolmogorov's contribution. For the third results I introduce two new meanings of `observational equivalence at every observation level'. Then I show the even stronger result of observational equivalence at every (and not just some) observation level between deterministic and indeterministic descriptions found in science. These results imply the following. Suppose one wants to find out whether a phenomenon is best modeled as deterministic or indeterministic. Then one cannot appeal to differences in the probability distributions of deterministic and indeterministic descriptions found in science to argue that one of the descriptions is preferable because there is no such difference. Finally, I criticise the extant claims of philosophers and mathematicians on observational equivalence.

\newpage
\tableofcontents
\newpage
\section{Introduction}

Determinism and indeterminism and whether one can know that phenomena are governed by deterministic or indeterministic laws are crucial philosophical themes. Hence it is an important question whether there is observational equivalence between deterministic and indeterministic descriptions. However, this question has hardly been discussed. This paper contributes to filling this gap by \textit{presenting and philosophically assessing three types of results of increasing strength on the observational equivalence of deterministic and indeterministic descriptions}. When saying that a deterministic and an indeterministic description are observationally equivalent, I mean that the deterministic description, when observed, gives the same predictions as the indeterministic description. The deterministic and indeterministic descriptions of concern are continuous-time measure-theoretic deterministic systems and stochastic processes.

More specifically, the first results are about observational equivalence to abstract mathematical descriptions.  I present a method of constructing, given a deterministic system, an observationally equivalent stochastic process, and conversely. Compared to the first results, the second results are stronger because they show not only observational equivalence to \textit{abstract mathematical descriptions} but observational equivalence between deterministic and stochastic descriptions \textit{of the types antecedently found in science}. Compared to the second results, the third results are even stronger because they show that there is observational equivalence \textit{at every (not just at some) observation level} between deterministic and stochastic descriptions found in science. Notice that the increase in strength is very different for the move from the first to the second results (observational equivalence between descriptions found in science and not just to abstract mathematical descriptions) and for the move from the second to the third results (observational equivalence between descriptions in science at every, and not just at some, observation level).

I argue that a philosophical consequence of these results is the following.  Suppose one wants to find out whether a phenomenon is best modeled as deterministic or stochastic.
Then one might think of arguing that there is evidence for a deterministic or stochastic description by appealing to a general difference between the probability distributions of stochastic processes found in science and the probability distributions of (possibly fine-enough) observations of deterministic systems found in science. The second and third results of this paper show that these arguments are untenable because there is no such general difference between the respective probability distributions.

Finally, I criticise the previous philosophical discussion. The main philosophy literature on this topic is Suppes and de~Barros~(1996), Suppes~(1999) and Winnie~(1998). They claim that the significance of the third results is to provide a choice between deterministic and stochastic descriptions. I argue that while there is indeed such a choice, this is already shown by more basic results, and the third results show something stronger. Also, I criticise the claims of the mathematicians Ornstein~and Weiss~(1991) on observational equivalence.

As mentioned above, this paper is about \textit{continuous-time} measure-theoretic deterministic systems and stochastic processes (here the time parameter varies continually). There are also \textit{discrete-time} measure-theoretic deterministic systems and stochastic processes (here the time parameter varies in discrete steps).
Werndl (2009a) discusses results about discrete-time descriptions which answer the question (for discrete time) of whether deterministic and indeterministic descriptions can be observationally equivalent. One of the contributions of my paper is to answer this question for continuous-time deterministic systems and stochastic processes.  This is important because, first, the discrete-time results leave open the answers for continuous time, and the answers for continuous time do not follow automatically from the discrete-time results. Indeed, as I will explain, the results differ for continuous-time and discrete-time descriptions; also, the proofs for continuous-time descriptions are harder and involve different techniques. Second, continuous-time descriptions are more widespread in science than discrete-time descriptions. Hence the issue of observational equivalence is more pressing for continuous-time descriptions.

Furthermore, there are several novel contributions of this paper which are not discussed in Werndl (2009a) or the other extant literature. In particular, I show that there are results on observational equivalence which increase in strength in a certain sense (namely, the focus on descriptions found in science and on every observation level). For this I introduce two new notions of observational equivalence at every observation level, I derive results for these notions, and I discuss Kolmogorov's contribution to observational equivalence. Moreover, I assess the significance of these results and point at their philosophical consequences. Also, the focus on continuous-time enables me to criticise the philosophical reflections of mathematicians, namely Ornstein and Weiss (1991). And I criticise in detail the claims of philosophers such as Suppes and de Barros~(1996), Suppes~(1999) and Winnie~(1998) on observational equivalence.

This paper proceeds as follows. In Section~\ref{1} I introduce continuous-time measure-theoretic deterministic descriptions and stochastic processes. In Section~\ref{BIC} I present and assess the first results, in Section~\ref{AI1} the second results, and in Section~\ref{AI2} the third results on observational equivalence.  Finally, in Section~\ref{PM2} I criticise the previous philosophical discussion on observational equivalence by philosophers and mathematicians.

\section{Continuous-Time Deterministic Systems and Stochastic Processes}\label{1}

For what follows, I need to introduce a few basic measure-theoretic notions. Intuitively speaking, a probability space $(M,\Sigma_{M},\mu)$ consists of a set $M$, a set $\Sigma_{M}$ of subsets of $M$ to which a probability is assigned, called a sigma-algebra of $M$, and a probability measure $\mu$ which assigns a probability to subsets of $M$.
Formally, $(M,\Sigma_{M},\mu)$ is a \textit{probability space} if, and only if (iff), $M$ is a set; $\Sigma_{M}$ is a \textit{sigma-algebra} of $M$, i.e., a set of subsets of $M$ with (i) $\emptyset\in\Sigma_{M}$, (ii) for all $A\in\Sigma_{M}$, $M\setminus A\in\Sigma_{M}$, and (iii) for any $A_{n}\in \Sigma_{M},\,n\geq 1,\,\,\bigcup_{n}A_{n}\in\Sigma_{M}$; and $\mu$ is a \textit{probability measure}, i.e., a function $\mu:\Sigma_{M}\rightarrow [0,1]$ with (i) $\mu(A)\geq 0$ for all $A\in\Sigma_{M}$, (ii) $\mu(\emptyset)=0$ and $\mu(M)=1$, and (iii) $\sum_{n\geq 1}\mu(A_{n})=\mu(\bigcup_{n\geq 1}A_{n})$ for any $A_{n}\in\Sigma_{M}$. A pair $(M,\Sigma_{M})$ where $\Sigma_{M}$ is a sigma-algebra of $M$ is called a \textit{measurable space}. Finally, mathematically treatable functions in measure theory are called measurable functions; the functions encountered in science are all measurable. Technically, a function $T:M\rightarrow N$, where $(M,\Sigma_{M})$ and $(N,\Sigma_{N})$ are measurable spaces, is \textit{measurable} iff $T^{-1}(A)\in\Sigma_{M}$ for all $A\in\Sigma_{N}$.\footnote{For more details on measure theory, see Cornfeld et al.~(1982), Doob~(1953) and Petersen~(1983).}

\subsection{\textit{Deterministic Systems}}\label{Det}
This paper is about continuous-time measure-theoretic deterministic systems, in short deterministic systems, which are widespread in science.
The three main elements of a deterministic system are the \textit{phase space}, i.e., a set $M$ of all possible states, a probability measure $\mu$ which assigns a probability to regions of phase space\footnote{
The question of how to interpret this probability is not one of the main topics here. I just mention a popular interpretation. According to the \textit{time-average interpretation}, the probability of a set $A$ is the long-run fraction of the proportion of time the system spends in $A$ (cf.~Eckmann and Ruelle~1985; Werndl~2009b).}, and the
 \textit{evolution functions}, i.e., a family of functions $T_{t}:M\rightarrow M$ where $T_{t}(m)$ represents the state of the system after $t$ time units that started in initial state $m$. Formally:
\begin{definition}
A $\mathrm{deterministic}$ $\mathrm{system}$ is a quadruple $(M,\Sigma_{M},\mu,T_{t})$ where $(M,\Sigma_{M},\mu)$ is a probability space ($M$ is the phase space) and $T_{t}:M\rightarrow M$ (the evolution functions), $t\in\field{R}$, are measurable functions such that $T_{t_{1}+t_{2}}(m)=T_{t_{2}}(T_{t_{1}}(m))$ for all $m\in M$ and all $t_{1},t_{2}\in\field{R}$, and $T(t,m):=T_{t}(m)$ is measurable in $(t,m)$.\end{definition}
A solution is a possible path of the deterministic system. Formally, the \textit{solution} through $m$, $m\in M$, is the function $s_{m}:\field{R}\rightarrow M,\,\,s_{m}(t)=T(t,m)$. All deterministic systems are \emph{deterministic} according to the canonical definition: any solutions that agree at one instant of time agree at all times (Butterfield~2005).

I will often deal with measure-preserving
deterministic systems where, intuitively speaking, the probability of any region remains the same when the region is evolved.
\begin{definition}\label{MPDS}
A $\mathrm{measure}$-$\mathrm{preserving}$ $\mathrm{deterministic}$ $\mathrm{system}$  is a deterministic
system \linebreak[4] $(M,\Sigma_{M},\mu,T_{t})$ where $\mu$ is
\textit{invariant}, i.e.,\ $\mu(T_{t}(A))=\mu(A)$ for all $A\in\Sigma_{M}$ and $t\in\field{R}$.
\end{definition}

Measure-preserving deterministic systems are important in all the sciences, especially in physics (cf.~Eckmann and Ruelle~1985). All deterministic Hamiltonian systems and deterministic systems in statistical mechanics are measure-preserving, and their invariant probability measure is the normalised Lebesgue measure or a similar measure (Cornfeld et al.~1982, 4--10). A measure-preserving deterministic system is \emph{volume-preserving} iff the normalised Lebesgue measure -- the standard measure of the volume -- is the invariant probability measure. A measure-preserving system which fails to be volume-preserving is \emph{dissipative}. Dissipative systems can also often be modeled as measure-preserving systems, e.g., the Lorenz system, which was introduced by Lorenz (1963) to model weather phenomena (cf.~Luzatto et al.~2005).

When observing a deterministic system $(M,\Sigma_{M},\mu,T_{t})$, one observes a value $\Phi(m)$ dependent on the actual state $m$. Thus, technically, observations are \textit{observation functions}, i.e.,~measurable functions $\Phi:M\rightarrow M_{O}$ where $(M_{O},\Sigma_{M_{O}})$ is a measurable space (cf.~Ornstein and Weiss 1991, 16). An observation function $\Phi$ is \textit{finite-valued} iff it takes only  finitely many values $\{o_{1},\ldots,o_{k}\}$, $k\in\field{N}$, and every $o_{i}$ can be ``seen'', i.e., $\mu(\{m\in M\,|\,\Phi(m)=o_{i}\})>0$ for all $i,\,\,1\leq i\leq k$. A finite-valued observation function is \textit{nontrivial} iff $k\geq 2$. In practice observations are finite-valued.

The following deterministic system will accompany us. \\

\begin{figure}
\centering
\includegraphics{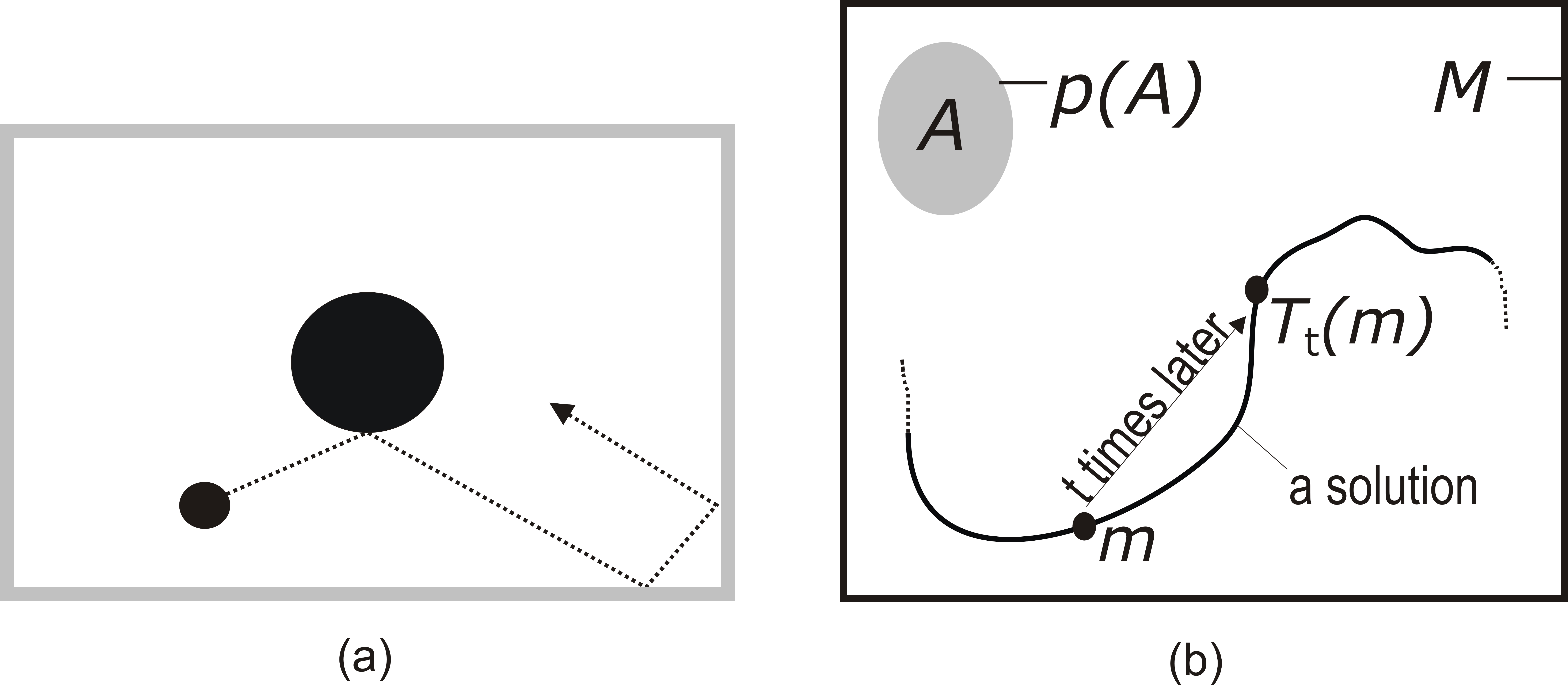}
\caption{A billiard system with a convex obstacle; (a) a specific state of the billiard (b) the mathematical representation of the billiard}\label{billiard}
\end{figure}

\noindent\textbf{Example~1: A billiard system with convex obstacles.}\\
A billiard system with convex obstacles is a system where a ball moves with constant speed on a rectangular table where there are a finite number of convex obstacles with a smooth boundary. It is assumed that there is no friction and that there are perfectly elastic collisions  (cf.~Ornstein and Galavotti~1974). Figure~1(a) shows a specific state of the billiard system, where the dashed line indicates the path the ball will take. Figure~1(b) shows the mathematical representation of the system. $M$ is the set of all possible positions and directions of the ball. Hence the state of the billiard shown in Figure~1(a) corresponds to exactly one point $m\in M$ shown in Figure~1(b). The normalised Lebesgue measure $\mu$ assigns the probability $\mu(A)$ to the event that the billiard is in one of the states represented by $A$ for any region $A$ in $\Sigma_{M}$ (intuitively speaking, $\Sigma_{M}$ is the set of all regions; formally, $\Sigma_{M}$ is the Lebesgue $\sigma$-algebra). $T_{t}(m)$, where $m=(q,p)$, gives the position and the direction after $t$ time units of the ball that starts out in initial position $q$ and initial direction $p$. And a solution is a possible path of the billiard over time.

\subsection{\textit{Stochastic Processes}}\label{SP}

The indeterministic descriptions of concern in this paper are continuous-time stochastic processes, in short stochastic processes. Stochastic processes are governed by probabilistic laws. They are the main indeterministic descriptions used in science.

A random variable $Z$ gives the outcome of a probabilistic experiment, where the distribution $P\{Z\in A\}$ tells one the probability that the outcome will be in $A$. Formally, a \textit{random variable} is a measurable function $Z:\Omega\rightarrow\bar{M}$ from a probability space $(\Omega,\Sigma_{\Omega},\nu)$ to a measurable space $(\bar{M},\Sigma_{\bar{M}})$. $P\{Z\in A\}=\nu(Z^{-1}(A))$ for all $A\in\Sigma_{\bar{M}}$ is the \textit{distribution} of $Z$.

A stochastic process $\{Z_{t};\,\,t\in\field{R}\}$ consists of a probabilistic experiment (i.e., a random variable) for each time $t$. So $Z_{t}(\omega)$ gives the outcome of the process at time $t$ (where $\omega$ represents a possible history in all its details). Formally:
\begin{definition}\label{stochproC}
A $\mathrm{stochastic}$ $\mathrm{process}$
$\{Z_{t};\,t\in\field{R}\}$ is a family of random
variables $Z_{t},\,\,t\in\field{R}$, which are defined on the same probability space $(\Omega,\Sigma_{\Omega},\nu)$ and take values
in the same measurable space $(\bar{M},\Sigma_{\bar{M}})$ such that $Z(t,\omega)=Z_{t}(\omega)$ is jointly measurable in $(t,\omega)$.
\end{definition}
$\bar{M}$ is the set of possible outcomes of the process, called the \textit{outcome space}. A \textit{realisation} is a possible path of the stochastic process. Formally, it is function $r_{\omega}:\field{R}\rightarrow\bar{M}$, $r_{\omega}(t)=Z(t,\omega)$, for $\omega\in\Omega$ arbitrary (cf.~Doob 1953,~4--46).

I will often deal with stationary processes, the probability distributions of which do not change with time. $\{Z_{t};\,t\in\field{R}\}$ is \textit{stationary} iff the distribution of $(Z_{t_{1}+h},\ldots, Z_{t_{n}+h})$ is the same as the one of $(Z_{t_{1}},\ldots, Z_{t_{n}})$ for all $t_{1},\ldots,t_{n}\in\field{R},\,n\in\field{N}$ and $h\in\field{R}$.

When observing a stochastic process at time $t$, one observes a value $\Gamma(Z_{t})$ dependent on the outcome $Z_{t}$. Hence, technically, observations are modeled by \textit{observation functions}, i.e., measurable functions $\Gamma:\bar{M}\rightarrow\bar{M}_{O}$, where $(\bar{M}_{O},\Sigma_{\bar{M}_{O}})$ is a measurable space. Clearly, observing the stochastic process $\{Z_{t};\,\,t\in\field{R}\}$ with $\Gamma$ yields the process $\{\Gamma(Z_{t});\,t\in\field{R}\}$.

The following stochastic processes will accompany us.\\

\begin{figure}
\centering
\includegraphics{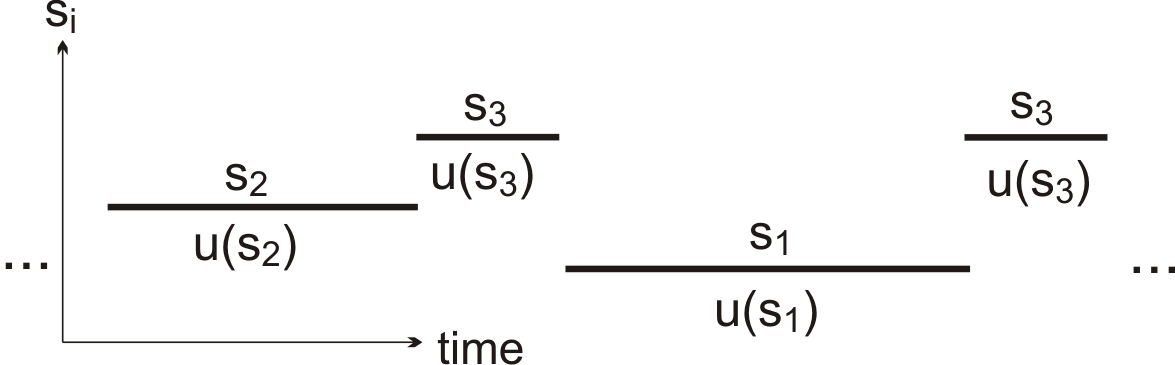}
\caption{A realisation of a semi-Markov process}\label{SMP}
\end{figure}

\noindent\textbf{Example~2: Semi-Markov processes.}\\
A semi-Markov process has finitely many possible outcomes $s_{i}$; it takes the outcome $s_{i}$ for a time $u(s_{i})$, and which outcome follows $s_{i}$ depends only on $s_{i}$ and no other past outcomes. Figure~2 shows a possible realisation of a semi-Markov process. For a formal definition, see Subsection~\ref{A1}. Semi-Markov processes are widespread in the sciences, from physics and biology to the social sciences. They are particularly important in queueing theory (cf.~Janssen and Limnios~1999).\\

\noindent\textbf{Example~3: $n$-step semi-Markov processes.}\\
$n$-step semi-Markov processes generalise semi-Markov processes.  An $n$-step semi-Markov process, $n\in\field{N}$, has finitely many possible outcomes $s_{i}$; it takes the outcome $s_{i}$ for a time $u(s_{i})$, and which outcome follows $s_{i}$ depends only on the past $n$ outcomes (hence semi-Markov processes are $1$-step semi-Markov processes). For a formal definition, see Subsection~\ref{A1}. Again, $n$-step semi-Markov processes are widespread in science (cf.~Janssen and Limnios~1999).\\

\noindent A final comment: the descriptions introduced in this section are classical. What results hold for quantum-mechanical descriptions and how similarities and differences in the results bear on the relation between classical and quantum mechanics are interesting questions; but they are beyond the scope of this paper.

\section{Observational Equivalence: Results I}\label{BIC}
This section is about observational equivalence to  mathematical descriptions. I show how, when starting with deterministic systems, one finds observationally equivalent stochastic processes (Subsection~\ref{DRS}), and how, when starting with  stochastic processes, one finds observationally equivalent deterministic systems (Subsection~\ref{SRD}).

I speak of observational equivalence if \textit{the deterministic system, when observed, and the stochastic process give the same predictions}. More specifically, deterministic systems are endowed with a probability measure. Hence when observing a deterministic system, the resulting predictions are the probability distributions over sequences of possible observations. The predictions obtained from a stochastic process are the probability distributions over the realisations. Consequently, a deterministic system, when observed, and a stochastic process give the same predictions iff (i) the possible outcomes of the process and the possible observed values of the system coincide, and (ii) the realisations of the process and the solutions of the system coarse-grained by the observation function have the same probability distribution.

\subsection{\textit{Starting With Deterministic Systems}}\label{DRS}
Assume $(M,\Sigma_{M},\mu,T_{t})$ is observed with $\Phi:M\rightarrow M_{O}$. Then $\{\Phi(T_{t});\,t\in\field{R}\}$ is a stochastic process, which arises by applying $\Phi$ to the deterministic system. Thus the possible outcomes of the process and the possible observed values of the deterministic system coincide, and the probability distributions over the realisations of the process and over the solutions of the system coarse-grained by the observation function are the same. Therefore, \textit{$(M,\Sigma_{M},\mu,T_{t})$ observed with $\Phi$ is observationally equivalent to $\{\Phi(T_{t});\,t\in\field{R}\}$}. \\

\begin{figure}
\centering
\includegraphics{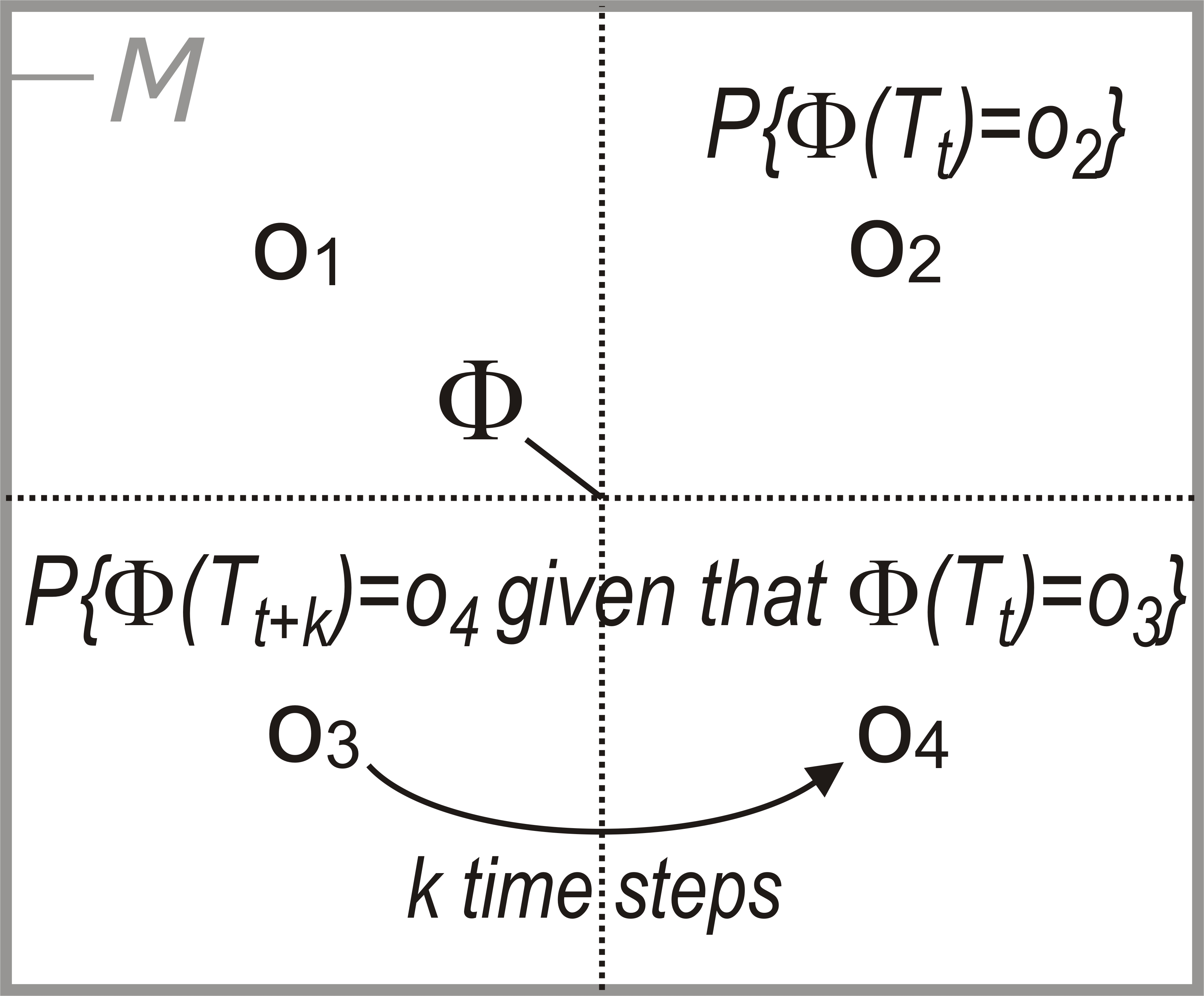}
\caption{The billiard observed with $\Phi$ is observationally equivalent to $\{\Phi(T_{t});\,t\in\field{R}\}$.}\label{todet}
\end{figure}

Let me illustrate this with the billiard system with convex obstacles. Figure~3 shows the phase space $M$ of the billiard system (cf.~Figure~1(b)) and an observation function $\Phi$ with values $o_{1},o_{2},o_{3},o_{4}$: for all states $m$ in the top left box the observed value is $o_{1}$, etc. Because a probability measure $\mu$ is defined on $M$, one obtains probabilities such as $P\{\Phi(T_{t})=o_{2}\}$ and $P\{\Phi(T_{t+k})=o_{4}$ given that $\Phi(T_{t})=o_{3}\}$. $\{\Phi(T_{t});\,\,t\in\field{R}\}$ has the outcomes $o_{1},o_{2},o_{3},o_{4}$ and its probability distributions are obtained by observing the billiard system with $\Phi$. Thus the billiard system observed with $\Phi$ is observationally equivalent to $\{\Phi(T_{t});\,t\in\field{R}\}$.

However, we want to know whether $\{\Phi(T_{t});\,t\in\field{R}\}$ is nontrivial. For let $\Phi(m)=m$ be the identity function. $\{\Phi(T_{t}(m));\,\,t\in\field{R}\}=\{T_{t}(m);\,\,t\in\field{R}\}$ is formally a stochastic process which is observationally equivalent to the deterministic system. Yet it is just the original deterministic system. This result is the formalisation of the old idea that a deterministic system is the special case of a stochastic process where all probabilities are $0$ or $1$ (cf.~Butterfield~2005). It illustrates that if one wants to arrive at nontrivial stochastic processes, one has to apply  observation functions which are many-to-one.

Yet \textit{$\{\Phi(T_{t});\,t\in\field{R}\}$ is often nontrivial}.  Several results show this. To my knowledge, the following theorem has never been discussed before. This theorem characterises a class of deterministic systems as systems where, regardless which finite-valued observation function $\Phi$ is applied, $\{\Phi(T_{t});\,t\in\field{R}\}$ is nontrivial in the following sense: for any $k\in\field{R}^{+}$ there are outcomes $o_{i},o_{j}\in M_{O}$  such that the probability of moving in $k$ time steps from $o_{i}$ to $o_{j}$ is between $0$ and $1$. This result is strong because there are nontrivial probability distributions for \textit{any} finite-valued observation function and \textit{all} time steps.

\begin{theorem}\label{epC}
Iff for a measure-preserving deterministic system $(M,\Sigma_{M},\mu,T_{t})$ there does not exist a $n\in\field{R}^{+}$ and a $C\in\Sigma_{M}$, $0<\mu(C)<1,$ such that, except for a set of measure zero\footnote{That is, except for $Q$ with $\mu(Q)=0$.}, $T_{n}(C)=C$, then the following holds: for every nontrivial finite-valued observation function $\Phi:M\rightarrow M_{O}$, every $k\in\field{R}^{+}$ and $\{Z_{t};\,t\in\field{R}\}:=\{\Phi(T_{t});\,t\in\field{R}\}$ there are $o_{i},o_{j}\in M_{O}$ with $0<P\{Z_{t+k}\!=\!o_{j}\,|\,Z_{t}\!=\!o_{i}\}<1$.
\end{theorem} For a proof, see Subsection~\ref{A2}.\footnote{If $(M,\Sigma_{M},\mu,T_{t})$ is measure-preserving, $\{\Phi(T_{t});\,\,t\in\field{R}\}$ is stationary (the proof in Werndl~2009a,~Section~3.3, carries over to continuous time).}

The assumption of Theorem~1 is equivalent to weak mixing (Hopf 1932). This indicates a difference between the results for discrete and continuous time. As shown in Werndl (2009a), for discrete time the condition needed to arrive at the results that, regardless which finite-valued observation function and time step is considered, one always obtains nontrivial probability distributions is \textit{weaker} than weak mixing. In particular, Werndl (2009a,~Section 4.2.2) shows that this result holds even for a class of \textit{stable} deterministic systems (i.e., where solutions which start closely stay close all the time). This is not the case for continuous time because if a deterministic system is weak mixing, there is a sense in which solutions which start closely eventually separate (cf.~Werndl~2009a). Yet Theorem~1 is just one of the possible theorems one can prove. Dependent on one's purpose, one can  weaken (or strengthen) the conditions to prove weaker (or stronger\footnote{For instance, a stronger result is that for Kolmogorov systems and any finite-valued observation function the following holds: even if you know the entire infinite history of the process, you will not be able to predict with certainty which outcome follows next (Uffink~2006,~1012--1014).}) results on observational equivalence.
And weaker results than Theorem 1 show that one obtains observational equivalence to nontrivial stochastic processes also for stable deterministic systems and thus also for non-chaotic systems (\textit{chaotic} systems show deterministic yet unstable behaviour). For instance, consider the deterministic system of a rotation on a circle, which is stable and non-chaotic; $M=[0,1)$ represents the circle of unit radius, $\Sigma_{M}$ is the Lebesgue $\sigma$-algebra, $\mu$ the Lebesgue measure and $T_{t}(m)=\alpha tm\,(\textnormal{mod}\,1)$, $\alpha\in\field{R}$. It is not hard to see that for any finite-valued observation function $\Phi$ the stochastic process $\{Z_{t};\,t\in\field{R}\}=\{\Phi(T_{t});\,t\in\field{R}\}$ is nontrivial in the sense that for almost all (but not all) time steps there are nontrivial probability distributions. Technically:
for any finite-valued $\Phi$ and all, except for a set of measure zero, $k\in\field{R}^{+}$, there are $o_{i},o_{j}\in M_{O}$ with $0<P\{Z_{t+k}\!=\!o_{j}\,|\,Z_{t}\!=\!o_{i}\}<1$, where $\{Z_{t}\}=\{\Phi(T_{t});\,t\in\field{R}\}$.

Measure-preserving systems are typically weakly mixing (Halmos~1944, 1949). Consequently, Theorem~1 typically applies to measure-preserving deterministic systems. Yet this does not tell us whether Theorem~1 applies to the deterministic systems of physical relevance; systems of physical relevance are only a small class of all measure-preserving  systems. But Theorem~\ref{epC} applies to several physically relevant deterministic systems. It applies to systems in Newtonian mechanics: for instance, first, to billiard systems with convex obstacles (Example~1) (Ornstein and Galavotti 1974); second, to many hard-ball systems, which describe the motion of a number of hard balls undergoing elastic reflections at the boundary and collisions amongst each other and are important in statistical mechanics because they model gases; e.g., to two hard balls in a box and the motion of $N$ hard balls on a torus for almost all values $(m_{1},\ldots,m_{N},r)$, where $m_{i}$ is the mass of the $i$-th ball and $r$ is the radius of the balls,
$N\geq 2$ (Berkovitz et al.~2006, 679--680; Sim\'{a}nyi~2003); third, to geodesic flows of negative curvature, i.e., frictionless motion of a particle moving with unit speed on a compact manifold with everywhere negative curvature (Ornstein and Weiss~1991, Section~4). Also, Theorem~\ref{epC} applies to dissipative systems, such as Lorenz-type systems, which model weather dynamics and waterwheels (Lorenz~1963; Luzzatto et al.~2005; Strogatz~1994). Because proofs are extremely hard, it is often only conjectured that deterministic systems satisfy Theorem~\ref{epC}, e.g., for any number of hard balls in a box (Berkovitz et al.~2006, 679--680). Note that the deterministic systems listed in this paragraph are chaotic.

There are also many systems to which Theorem~\ref{epC} does not apply. For instance, according to the KAM theorem, the phase space of integrable Hamiltonian systems which are perturbed by a small non-integrable perturbation breaks up into stable regions and regions with unstable behavior. Because of the break up into regions, Theorem~1 does not apply to these systems (cf.~Berkovitz et al.~2006, Section~4). But even if Theorem~\ref{epC} does not apply to the whole deterministic system, it might well apply to the motion restricted to some regions $A$ of phase space (and this is conjectured for KAM-type systems). Then it follows from Theorem~\ref{epC} that all finite-valued observations which discriminate between values in $A$ lead to nontrivial stochastic processes.

\subsection{\textit{Starting With Stochastic Processes}}\label{SRD}
The following idea of how, given stochastic processes, one finds observationally equivalent deterministic systems will be needed for what follows. It is well known in mathematics and in philosophy of physics (Butterfield~2005; Doob~1953, 6--7; for discrete time it is discussed in Werndl~2009a, Section~3.2).

The underlying idea is that one constructs a deterministic system whose phase space is the set of all possible realisations $m(\tau)$, whose evolution functions $T_{t}$ shift the realisation $t$ time steps to the left, and whose observation function $\Phi_{0}$ gives the value of the realisation at time $0$. Because $T_{t}$ shifts the realisation $m(\tau)$ $t$ time steps to the left, observing $T_{t}(m(\tau))$ with $\Phi_{0}$ gives the value of the realisation $m(\tau)$ at time $t$. Hence successive observations return the values of the realisation $m(\tau)$ at different times and so simply yield the outcomes of the stochastic process.

Formally, for $\{Z_{t};\,t\in\field{R}\}$ with outcome space $\bar{M}$, let $M$ be the set of all possible realisations, i.e., functions $m(\tau)$ from $\field{R}$ to $\bar{M}$. Let $\Sigma_{M}$ be the $\sigma$-algebra generated\footnote{The $\sigma$-algebra generated by $E$ is the smallest $\sigma$-algebra containing $E$; that is, the $\sigma$-algebra $\bigcap_{\sigma\textnormal{-algebras}\,\Sigma,\,E\subseteq \Sigma}\Sigma$.}
by the cylinder-sets
\begin{equation}\label{cylinder}
C^{A_{1}...A_{n}}_{i_{1}...i_{n}}\!\!=\!\!\{m\!\in\! M\,|\,m(i_{1})\!\!\in\!\! A_{1},...,m(i_{n})\!\!\in\!\! A_{n}, A_{j}\!\in\! \Sigma_{\bar{M}},i_{j}\!\in\!\field{R},\,i_{1}\!\!<...<\!\!i_{n},1\!\leq j\!\leq n\}.
\end{equation}
Let $\mu$ be the measure on $M$ determined by the probability distributions of  $\{Z_{t};\,\,t\in\field{R}\}$. That is, $\mu$ is the unique probability measure arising by assigning to each   $C^{A_{1}...A_{n}}_{i_{1}...i_{n}}$ the probability $P\{Z_{i_{1}}\in
A_{1},\ldots,Z_{i_{n}}\in A_{n}\}$. The evolution functions shift a realisation $t$ times to the left, i.e., $T_{t}(m(\tau))=m(\tau+t)$. $(M,\Sigma_{M},\mu,T_{t})$ is a deterministic system called the \textit{deterministic representation} of $\{Z_{t};\,t\in\field{R}\}$ (cf.~Doob 1953,~621--622). Finally, assume that the observation function gives the value of the realisation at time zero, i.e., $\Phi_{0}(m(\tau))=m(0)$.

Because one only observes the 0-th coordinate, the possible outcomes of a stochastic process and the possible observed values of its deterministic representation coincide. Moreover, because $\mu$ is defined by the stochastic process, the probability distribution over the realisations is the same as the one over the sequences of observed values of its deterministic representation. Therefore, \textit{a stochastic process is observationally equivalent to its
deterministic representation observed with $\Phi_{0}$. Consequently, for every process there is at least one observationally equivalent deterministic system.}\footnote{The deterministic representation of any stationary process is measure-preserving (the proof in Werndl~2009a, Section~3.3, carries over to continuous time).}
For instance, consider a semi-Markov process $\{Z_{t};\,\,t\in\field{R}\}$ (Example~2) and its deterministic representation $(M,\Sigma_{M},\mu,T_{t})$. Here $M$ is the set of all possible realisations of the semi-Markov process, $\mu$ gives the probability distribution over the realisations, $T_{t}$ shifts a realisation $t$ time steps to the left, and $\Phi_{0}$ returns the value of the realisation at time $0$.

\textit{Philosophically speaking, the deterministic representation is a cheat}. For it resembles the following argument: Fred was fatally wounded at noon; therefore he died later. The predicate `being fatally wounded a noon' implies that Fred will die later. Hence the argument is trivially valid. Similarly, the states of the deterministic representation are constructed such that they \textit{encode} the future and past evolution of the process. So there is the question whether deterministic systems which do not involve a cheat can be observationally equivalent to a given stochastic process. It will become clear later that for some stochastic processes the answer is positive. It is unknown whether the answer is positive for every stochastic process. For a formalisation of the notion of observational equivalence underlying Section~\ref{BIC}, see Werndl~2009a (the definition carries over to continuous-time).

\section{Observational Equivalence: Results II}\label{AI1}
The previous section only showed that in certain cases there is observational equivalence to abstract mathematical descriptions. Furthermore, the deterministic representation -- one of these abstract mathematical descriptions -- is not used by scientists. This raises the question of whether there can be observational equivalence not only to \textit{abstract mathematical descriptions} but whether a \textit{stronger} result can hold: namely, that there is observational equivalence between stochastic processes and deterministic systems \textit{of the types used in scientific theorising }(in short, \textit{systems and processes in science}).

If the answer is negative, one could divide the probability distributions found in science into two groups: the ones deriving from observations of deterministic systems in science, and the ones deriving from stochastic processes in science and their observations.\footnote{By `probability distributions' I mean here nontrivial probability distributions (not all probabilities $0$ or $1$). Otherwise such a division could never exist.
For consider a process composed of two semi-Markov processes (Example~2), and an observation function which tells one whether the outcome is among those of the first or the second process. Then the observed probability distribution is trivial.} Suppose one wants to find out whether a phenomenon is best described as deterministic or stochastic. Then one might argue the following: if the observed probability distributions are of the type of deterministic systems in science, the evidence speaks for a deterministic description. Conversely, if the probability distributions are of the type of stochastic processes in science, this provides evidence for a stochastic description. The idea here is that the probability distributions characteristically found for deterministic (or stochastic) descriptions in science provide inductive evidence for a deterministic (or stochastic) description. There is the question whether such an argument is convincing. Clearly, it can only work if there is indeed a division between the respective probability distributions. So let us ask: \textit{is there a deterministic system in science and is there an observation function such that the deterministic system, when observed, yields a stochastic process in science?}

\subsection{\textit{Kolmogorov's Conjecture}}
Kolmogorov found it hard to imagine a positive answer to this question.  More specifically, it has been nearly forgotten (and has neither been subject to a systematic historical investigation nor been discussed philosophically) that Kolmogorov conjectured the following. \textit{Suppose an observed deterministic system yields a stochastic process in science; then the observed deterministic system produces positive information. Contrary to this, any arbitrary observation of a deterministic system in science never produces positive information.} Employing ideas in information theory, the Kolmogorov-Sinai entropy (KS-entropy) was introduced to capture the highest amount of information, or equivalently the highest amount of uncertainty, that an observation of a deterministic system could produce; where a positive KS-entropy indicates that there are observations which produce positive information. Kolmogorov expected that the KS-entropy would accomplish the separation of deterministic systems in science from the deterministic systems producing stochastic processes in science: the former have zero entropy, the latter have positive entropy.

So it was a big surprise when, from the 1960s onwards, it was found that also many deterministic systems in science have positive KS-entropy and thus have observations which produce positive information. Namely, all the deterministic systems listed in Subsection~\ref{DRS} as systems of physical relevance to which Theorem~1 applies have positive KS-entropy; in particular, billiard systems with convex obstacles (Example~1), many hard-ball systems, geodesic flows of constant negative curvature and Lorenz-type systems. For many systems in science, such as for KAM-type systems, it is conjectured that they have positive KS-entropy. \textit{Hence Kolmogorov's attempt to separate deterministic systems in science from deterministic systems which produce stochastic processes in science failed} (Radunskaya~1992, chapter 1; Sinai~1989, 835--837).

\subsection{\textit{Deterministic Systems in Science Which Are Observationally Equivalent to Stochastic Processes in Science}}\label{eaglewing}
So is there a deterministic system in science and an observation such that the observed deterministic system yields a stochastic process in science? To answer this, I need to introduce the definitions of a continuous Bernoulli system and of isomorphism.

Intuitively speaking, Bernoulli systems are strongly chaotic.
\begin{definition}\label{cbSS}
The measure-preserving deterministic system $(M,\Sigma_{M},\mu,T_{t})$ is a $\mathrm{continuous}$ $\mathrm{Bernoulli}$ $\mathrm{system}$ iff for all $t\in\field{R}\setminus\{0\}$ the discrete measure-preserving deterministic system  $(M,\Sigma_{M},\mu,T_{t})$ is a discrete Bernoulli system.
\end{definition}
For the discussion in the main text, the details of this definition are not important.\footnote{A discrete Bernoulli system is defined as follows. Recall Definition~\ref{DS} of a discrete measure-preserving deterministic system. A doubly-infinite sequence of independent rolls of an $N$-sided die where the probability of obtaining $s_{k}$ is $p_{k}$, $1\leq k\leq N$, $\sum_{k=1}^{N}p_{k}\!=\!1$, is a Bernoulli process.
Let $M$ be the set of all bi-infinite sequences $m=(\ldots
m(-1),m(0),m(1)\ldots)$ with
$m(i)\in\bar{M}=\{s_{1},\ldots,s_{N}\}$. Let $\Sigma_{M}$ be the $\sigma$-algebra generated by the cylinder-sets $C^{\{s_{l_{1}}\}...\{s_{l_{n}}\}}_{i_{1}\ldots i_{n}}$ as defined in (\ref{cylinder}) with $i_{j}\in\field{Z},s_{l_{j}}\in\bar{M}$. These sets have probability
$\bar{\mu}(C^{\{s_{l_{1}}\}...\{s_{l_{n}}\}}_{i_{1}\ldots i_{n}})=p_{s_{l_{1}}}...p_{s_{l_{n}}}$. Let $\mu$ the unique measure determined by $\bar{\mu}$ on $\Sigma_{M}$. Let
$T:M\rightarrow M,\,\,T((\ldots m(i)\ldots))=(\ldots m(i+1)\ldots)$. $(M,\Sigma_{M},\mu,T)$ is called a Bernoulli shift. Finally, a discrete measure-preserving deterministic system is a \textit{Bernoulli system} iff it is isomorphic to some Bernoulli shift (isomorphism is exactly defined as for continuous time -- see Definition~\ref{isomorphic}).}
 Here it just matters that several deterministic systems in science are continuous Bernoulli system; namely, all the systems listed in the previous subsection as systems with positive KS-entropy, in particular, billiard systems with convex obstacles (Example~1), many hard-ball systems, geodesic flows of constant negative curvature and Lorenz-type systems. Again, for several systems, e.g., for the motion on unstable regions of KAM-type systems, it is conjectured that they are continuous Bernoulli systems.

Isomorphic deterministic systems are probabilistically equivalent, i.e., their states can be put into one-to-one correspondence (via a function $\phi$) such that the corresponding solutions have the same probability distributions.
\begin{definition}\label{isomorphic}
The measure-preserving deterministic systems  $(M_{1},\Sigma_{M_{1}},\mu_{1},T^{1}_{t})$ and \linebreak[4]
$(M_{2},\Sigma_{M_{2}},\mu_{2},T^{2}_{t})$ are $\mathrm{isomorphic}$ iff there are  $\hat{M}_{i}\subseteq M_{i}$ with
$\mu_{i}(M_{i}\setminus \hat{M}_{i})=0$ and
$T^{i}_{t}\hat{M}_{i}\subseteq\hat{M}_{i}$ for all $t\,\,(i=1,2$), and there is a
bijection $\phi:\hat{M}_{1}\!\rightarrow\!\hat{M}_{2}$ such that
(i) $\phi(A)\!\in\!\Sigma_{M_{2}}$ for all
$A\!\in\!\Sigma_{M_{1}},A\subseteq \hat{M}_{1}$, and
$\phi^{-1}(B)\in\Sigma_{M_{1}}$ for all
$B\in\Sigma_{M_{2}},B\subseteq\hat{M}_{2}$; (ii)
$\mu_{2}(\phi(A))=\mu_{1}(A)$ for all
$A\in\Sigma_{M_{1}},\,A\subseteq\hat{M}_{1}$; (iii) $\phi(T^{1}_{t}(m))=T^{2}_{t}(\phi(m))$ for all $m\in\hat{M}_{1}, t\in\field{R}$.
\end{definition}
Now assume that $(M,\Sigma_{M},\mu,T_{t})$ is isomorphic (via $\phi:\hat{M}\rightarrow\hat{M_{2}}$) to the deterministic
representation $(M_{2},\Sigma_{M_{2}},\mu_{2},T^{2}_{t})$ of the stochastic process $\{Z_{t};\,t\in\field{R}\}$. This means that there is a one-to-one correspondence between the states of the system and the realisations of the process. Recall that observing $(M_{2},\Sigma_{M_{2}},\mu_{2},T^{2}_{t})$ with $\Phi_{0}$ yields $\{Z_{t};\,t\in\field{R}\}$ (where $\Phi_{0}(m_{2}(\tau))=m_{2}(0)$). Consequently, $(M,\Sigma_{M},\mu,T_{t})$ observed with the observation function $\Phi:=\Phi_{0}(\phi(m))$ is observationally equivalent to $\{Z_{t};\,t\in\field{R}\}$\footnote{It does not matter how $\Phi_{0}(\phi(m))$ is defined for $m\in M\setminus\hat{M}$.}

An important application of this principle is as follows. A deep result shows that up to a scaling of time~$t$ any continuous Bernoulli systems are isomorphic. That is, given continuous Bernoulli systems $(M,\Sigma_{M},\mu,T_{t})$ and $(M_{2},\Sigma_{M_{2}},\mu_{2},T_{t}^{2})$ there is a $c>0$ such that $(M,\Sigma_{M},\mu,T_{t})$ and $(M_{2},\Sigma_{M_{2}},\mu_{2},T_{ct}^{2})$ are isomorphic (Ornstein~1974). Now recall semi-Markov processes (Example~2).
Ornstein~(1970) proved that the deterministic representation $(M_{2},\Sigma_{M_{2}},\mu_{2},T^{2}_{t})$ of any semi-Markov process  is a continuous Bernoulli system. Clearly, for any $c\in\field{R}^{+}$, $(M_{2},\Sigma_{M_{2}},\mu_{2},T^{2}_{ct})$ is still the deterministic representation of a semi-Markov process. Therefore, given any continuous Bernoulli system in science $(M,\Sigma_{M},\mu,T_{t})$, there is a semi-Markov process whose deterministic representation $(M_{2},\Sigma_{M_{2}},\mu_{2},T^{2}_{ct})$ is isomorphic (via $\phi$) to $(M,\Sigma_{M},\mu,T_{t})$, implying that $(M,\Sigma_{M},\mu,T_{t})$ observed with $\Phi:=\Phi_{0}(\phi(m))$ produces the semi-Markov process. For instance, billiards with convex obstacles (Example~1) are continuous Bernoulli systems. Hence there is an observation function $\Phi$ which, applied to the billiard system, yields a semi-Markov processes.

The converse is also true: given any semi-Markov process, there is a Bernoulli system in science which is observationally equivalent to the semi-Markov process. For given any Bernoulli system in science $(M,\Sigma_{M},\mu,T_{t})$, $(M,\Sigma_{M},\mu,T_{ct})$ is also a system in science (e.g., for billiard systems this means that the constant speed of the ball is changed from $s$ to $cs$). Thus given any semi-Markov process there is a deterministic system in science $(M,\Sigma_{M},\mu,T_{ct})$ which is isomorphic to the deterministic representation $(M_{2},\Sigma_{M_{2}},\mu_{2},T^{2}_{t})$ of the semi-Markov process, implying that $(M,\Sigma_{M},\mu,T_{ct})$ observed with $\Phi:=\Phi_{0}(\phi(m))$ yields the semi-Markov process. For instance, given any semi-Markov process (Example~2), one can change the constant speed of a billiard system with convex obstacles such that an observation of the billiard system yields the semi-Markov process.

To sum up, deterministic systems in science (namely, certain Bernoulli systems) are observationally equivalent to semi-Markov processes (which are widespread in science).\footnote{$n$-step semi-Markov processes are  continuous Bernoulli systems too (Ornstein~1974). Therefore, an analogous argument shows that Bernoulli systems in science are observationally equivalent to $n$-step semi-Markov processes.}  For discrete time Werndl (2009a, Section~4) showed that some deterministic systems in science are observationally equivalent to Bernoulli processes or discrete-time Markov processes (both are widespread in science). Yet this left open whether continuous-time deterministic systems in science can yield stochastic processes in science. Indeed, Werndl's results on Bernoulli or Markov processes do not carry over to continuous time:  Bernoulli processes are only defined for discrete time and there are no deterministic systems in science which yield continuous-time Markov processes. The latter holds because a deterministic system can only be observationally equivalent to a continuous-time Markov process if the KS-entropy of the Markov process is not higher than the KS-entropy of the system; but deterministic systems in science have finite KS-entropy, and continuous-time Markov processes have infinite KS-entropy (Feldman and Smorodinsky~1971; Ornstein and Weiss~1991,~19).

So I conclude that the answer to the question advanced at the beginning of this section is positive: \textit{there are deterministic systems in science and observation functions such that the deterministic systems, when observed, yield stochastic processes in science.}\footnote{My arguments allow any meaning of `deterministic systems deriving from scientific theories' that
excludes the deterministic representation but is wide enough to include some continuous Bernoulli systems.} So there is observational equivalence
not just to abstract mathematical descriptions (as shown in Section~\ref{BIC})\footnote{Also, compared to the results in Section~\ref{BIC}, the results here concern a narrower class of descriptions, namely only continuous Bernoulli systems, semi-Markov processes and $n$-step semi-Markov processes.} but  a \textit{stronger} form of observational equivalence between descriptions in science. Thus any argument that a phenomenon is best described as deterministic or stochastic (such as the one sketched at the beginning of this section) cannot work when it relies on the premise that there is a division between the probability distributions of deterministic and stochastic descriptions in science.

\section{Observational Equivalence: Results III}\label{AI2}
In the previous section we have seen that there are \textit{some} observations of deterministic systems in science which yield stochastic processes in science. Now let us ask whether an even \textit{stronger} result holds. Namely, \textit{can deterministic systems in science and stochastic processes in science be observationally equivalent at every observation level} (and not just for some observation)? Note that, compared to the previous section, the type of systems considered remain the same, namely deterministic and stochastic descriptions in science.

Intuitively, one might think that deterministic systems in science can  yield stochastic processes in science \textit{only} if \textit{specific coarse} observation functions are applied. If this is true, fine-enough observations of deterministic systems in science would yield  different probability distributions than stochastic processes in science. Then, if one wants to find out whether a phenomenon is best described as deterministic or stochastic, one might argue the following: if fine observations yield  probability distributions characteristically found for deterministic systems in science, the evidence favours a deterministic description. The idea here is that the existence of characteristic probability distributions for deterministic systems found in science gives us inductive evidence that a phenomenon is best described as deterministic. Such an argument is only tenable if fine-enough observations of deterministic systems indeed yield special probability distributions, i.e., if the question whether deterministic systems in science can be observationally equivalent at every observation level to stochastic processes in science has a negative answer. So let us focus on this question.

\subsection{\textit{The Meaning of Observational Equivalence at Every Observation Level}}\label{muede}
What does it mean that `stochastic processes of a certain type are observationally equivalent to a deterministic system \textit{at every  observation level}'? I will first introduce the standard notion and then two new notions. I focus on measure-preserving deterministic systems and, correspondingly, stationary stochastic processes because all examples will be of this kind.

\subsubsection*{The usual meaning based on $\varepsilon$-congruence}\label{UM} To introduce the standard notion, I start by explaining what it means for a deterministic system and a stochastic process to give the same predictions at an observation level $\varepsilon>0,\,\,\varepsilon\in\field{R}$. For sufficiently small $\varepsilon_{1}>0$ one cannot
distinguish states of the deterministic system which are less than
the distance $\varepsilon_{1}$ apart. Also, for sufficiently small $\varepsilon_{2}>0$ one will not be able to distinguish
differences in probabilities of less than $\varepsilon_{2}$. Assume that $\varepsilon<\min\{\varepsilon_{1},\varepsilon_{2}\}$.
Then a deterministic system and a stochastic process give the same predictions at level $\varepsilon$ iff the following shadowing result holds: the solutions can be put into one-to-one
correspondence with the realisations such that the state of the  system and the outcome of the process are at all time points within distance $\varepsilon$ except for a set of probability smaller than $\varepsilon$.

$\varepsilon$-congruence captures this idea mathematically (one assumes that a metric $d_{M}$ measures distances between states and that the possible outcomes of the process are in $M$; also, recall the deterministic representation as discussed in Subsection~\ref{SRD} and Definition~\ref{isomorphic} of isomorphism):
\begin{definition}\label{schatzi}
Let $(M,\!\Sigma_{M},\!\mu,\!T_{t})$ be a measure-preserving
deterministic system, where $(M,\!d_{M})$ is a metric space.
Let $(M_{2},\Sigma_{M_{2}},\mu_{2},T^{2}_{t})$ be the
deterministic representation of the stochastic process
$\{Z_{t};\,t\in\field{R}\}$ and let $\Phi_{0}:M_{2}\rightarrow M,\,\,\Phi_{0}(m_{2}(\tau))=m_{2}(0)$. $(M,\Sigma_{M},\mu,T_{t})$ is
$\varepsilon$-$\mathrm{congruent}$ to $\{Z_{t};\,t\in\field{R}\}$ iff
$(M,\Sigma_{M},\mu,T_{t})$ is isomorphic via $\phi:M\rightarrow M_{2}$ to
$(M_{2},\Sigma_{M_{2}},\mu_{2},T^{2}_{t})$ and
$d_{M}(m,\Phi_{0}(\phi(m)))\!\!<\!\!\varepsilon$ for all $m\!\!\in\!\!M$
except for a set of measure $\!<\!\varepsilon$.
\end{definition}

Now we generalise over $\varepsilon$ and arrive at a plausible meaning of the phrase that `stochastic processes of a certain type are observationally equivalent at every observation level to a deterministic system'. Namely: \textit{for every $\varepsilon>0$ there is a process of this type which is $\varepsilon$-congruent to the system} (cf.~Ornstein and Weiss 1991, 22--23; Suppes~1999).

The idea for new notions of observational equivalence at every observation level starts from the following thought.
$\varepsilon$-congruence does \textit{not} assume that the system is observed with an observation function: the \textit{actual} states of the deterministic system, \textit{not} the observed ones, are compared with the outcomes of the stochastic process. To arrive at a notion of observational equivalence, no observation functions are invoked, but it is asked whether the state of the system and the outcome of the process are less than $\varepsilon$ apart. (I should mention that if $(M,\Sigma_{M},\mu,T_{t})$ and $\{Z_{t};\,\,t\in\field{R}\}$ are $\varepsilon$-congruent, then $\{\Psi(T_{t});\,\,t\in\field{R}\}$, where $\Psi(m):=\Phi_{0}(\phi(m))$, is the process $\{Z_{t};\,\,t\in\field{R}\}$. Technically, $\Psi$ is an observation function but for $\varepsilon$-congruence it is \textit{not} interpreted in this way. Instead, the meaning of $\Psi$ is as follows: when it is applied to the system, the resulting process shadows the deterministic system.)

In many contexts, such as in Newtonian and statistical mechanics, observations are modeled by observation functions, and one would like to know what happens when specific observation functions are applied. Thus it would be desirable to have a notion of observational equivalence at every observation level based on observation functions. Results about other notions might be regarded as wanting because they leave unclear what happens when observation functions are applied. As explained, $\varepsilon$-congruence is not based on observation functions.  For this reason, I now introduce two other notions of observational equivalence at every observation level. Whether a notion is preferable that (i) tells one what happens when one cannot distinguish between states which are less than $\varepsilon$ apart (such as the notion based on $\varepsilon$-congruence) or (ii) tells one what stochastic processes are obtained if specific observation functions are applied (such as the notions introduced in the next two paragraphs) will depend on the modeling process and the phenomenon under consideration.

\subsubsection*{A new meaning based on strong $(\Phi,\varepsilon)$-simulation}
First, I have to explain what it means for a stochastic process and a deterministic system as observed with an observation function $\Phi$ to give the same predictions relative to accuracy $\varepsilon>0,\,\,\varepsilon\in\field{R}^{+}$ ($\varepsilon$ indicates that one cannot distinguish differences in probabilistic predictions of less than $\varepsilon$). Plausibly, this means that the possible observed values of the system and the possible outcomes of the process coincide, and that the probabilistic predictions of the observed deterministic system and of the stochastic process differ by less than $\varepsilon$.
Strong $(\Phi,\varepsilon)$-simulation captures this idea mathematically. I assume that $\Phi$ is finite-valued, as it is in practice.
\begin{definition}\label{Schatzi2}
$\{Z_{t};\,t\in\field{R}\}$ $\mathrm{strongly}$ $(\Phi,\varepsilon)$-$\mathrm{simulates}$ a measure-preserving deterministic system $(M,\Sigma_{M},\mu,T_{t})$, where $\Phi:M\rightarrow\bar{M}$ is a surjective finite-valued observation function, iff there is a surjective measurable function $\Psi:M\rightarrow\bar{M}$ such that (i) $Z_{t}=\Psi(T_{t})$ for all $t\in\field{R}$ and (ii) $\mu(\{m\in\!M\,|\,\Psi(m)\neq\Phi(m)\})<\varepsilon$.
\end{definition}

By generalising over $\Phi$ and $\varepsilon$ one obtains a plausible meaning of the notion that processes of a certain type are observationally equivalent at every observation level to a system. Namely: for every finite-valued $\Phi$ and every $\varepsilon$ there is a stochastic process of this type which strongly $(\Phi,\varepsilon)$-simulates the deterministic system. This notion is attractive because it tells us what probability distributions are obtained when any finite-valued observation function is applied. Yet, to my knowledge, it has not been discussed before.

\subsubsection*{A new meaning based on weak $(\Phi,\varepsilon)$-simulation}
If one asks whether the probability distributions obtained by applying $\Phi$ to the deterministic system could derive from a stochastic process of a certain kind, the notion of strong $(\Phi,\varepsilon)$-simulation is stronger than what is needed. Is suffices that the observed probability distributions could result from an \textit{observation} of a stochastic process of a certain kind (this will be of some relevance later -- see the end of Subsection~\ref{bermading}).
More specifically, if suffices to require the following. Given a deterministic system observed with $\Phi$, there is an observation $\Gamma$ of a stochastic process such that the following holds: the possible observed outcomes of the process are the possible observed values of the deterministic system, and  the probabilistic predictions of the process observed with $\Gamma$ and  the probabilistic predictions of the system observed with $\Phi$ differ by less than $\varepsilon$. Technically, weak $(\Phi,\varepsilon)$-simulation captures this idea.
\begin{definition}\label{Schatzi3}  $\{Z_{t};\,t\in\field{R}\}$ $\mathrm{weakly}$ $(\Phi,\varepsilon)$-$\mathrm{simulates}$
 a measure-preserving deterministic system $(M,\Sigma_{M},\mu,T_{t})$, where $\Phi:M\rightarrow\bar{M}$ is a surjective finite-valued observation function, iff there is a surjective measurable function $\Psi:M\rightarrow S$ and a surjective observation function $\Gamma:S\rightarrow\bar{M}$ such that (i) $\Gamma(Z_{t})=\Psi(T_{t})$ for all $t\in\field{R}$ and (ii) $\mu(\{m\in M\,|\,\Psi(m)\neq \Phi(m)\})<\varepsilon$.
\end{definition}

Clearly, if a stochastic process strongly $(\Phi,\varepsilon)$-simulates a system, it also weakly $(\Phi,\varepsilon)$-simulates it (choose $\Gamma(s)=s$). The converse is generally not true. By generalising over $\Phi$ and $\varepsilon$, one obtains a plausible meaning of the notion that stochastic processes of a certain type are observationally equivalent at every observation level to a deterministic system. Namely: for every finite-valued $\Phi$ and every $\varepsilon$ there is a process of this type which weakly $(\Phi,\varepsilon)$-simulates the system. To my knowledge, this notion has never been discussed.

According to all three notions, at every observation level, the data could derive from the deterministic system or a stochastic process of a certain type. So let us see what results obtain for the three notions.

\subsection{\textit{Stochastic Processes in Science Which Are Observationally Equivalent at Every Observation Level to Deterministic Systems in Science}}\label{bermading}
The next three theorems show that the following holds for our three notions: continuous Bernoulli systems (Definition~\ref{cbSS}), including several systems in science, are observationally equivalent at every observation level to semi-Markov processes (notion one and three) or $n$-step semi-Markov processes (notion two), which are widespread in science.
\begin{theorem}\label{T5}
Let $(M,\Sigma_{M},\mu,T_{t})$ be a continuous Bernoulli system where the metric space $(M,d_{M})$ is separable\footnote{$(M,d_{M})$ is \textit{separable}
iff there is a countable set
$\tilde{M}=\{m_{n}\,|n\in\field{N}\}$ with $m_{n}\in M$ such that
every nonempty open subset of $M$ contains at least one element of
$\tilde{M}$.} and $\Sigma_{M}$ contains all open sets of $(M,d_{M})$. Then for every $\varepsilon>0$, $(M,\Sigma_{M},\mu,T_{t})$ is $\varepsilon$-congruent to a semi-Markov process.
\end{theorem}
For a proof see Ornstein and Weiss (1991,~93--94). The assumptions of this theorem are fulfilled by all continuous Bernoulli systems in science. For discrete time Werndl (2009a,~Section~4) showed that several deterministic systems in science are $\varepsilon$-congruent for all $\varepsilon>0$ to discrete-time Markov processes (which are widespread in science). However, this left open whether continuous-time deterministic systems in science can be $\varepsilon$-congruent for all $\varepsilon>0$ to continuous-time stochastic processes in science. Indeed, Werndl's~(2009a) results on Markov processes do not carry over to continuous time because, as explained at the end of Subsection~\ref{eaglewing}, deterministic systems in science cannot yield continuous-time Markov processes.

\begin{theorem}\label{T6}
Let $(M,\Sigma_{M},\mu,T_{t})$ be a continuous Bernoulli system. Then for every finite-valued observation function $\Phi$ and every $\varepsilon>0$ there is an $n$ such that an $n$-step semi-Markov process strongly $(\Phi,\varepsilon)$-simulates $(M,\Sigma_{M},\mu,T_{t})$.
\end{theorem}
For a proof see Ornstein and Weiss (1991,~94--95).

\begin{theorem}\label{semiMarkov1}
Let $(M,\Sigma_{M},\mu,T_{t})$ be a continuous Bernoulli system. Then for every finite-valued observation function $\Phi$ and every $\varepsilon>0$ a semi-Markov process $\{Z_{t},\,t\in\field{R}\}$ weakly  $(\Phi,\varepsilon)$-simulates $(M,\Sigma_{M},\mu,T_{t})$.
\end{theorem}
For a proof, see Subsection~\ref{P8}.

For instance, consider a billiard system with convex obstacles (Example~1). According to Theorem~2, the solutions of such a billiard system are shadowed by the realisations of a semi-Markov process for any accuracy $\varepsilon$. According to Theorem~3, such a billiard system, observed with any finite-valued observation function, is observationally equivalent to an $n$-step semi-Markov process (disregarding differences in probabilistic predictions of less than $\varepsilon$). According to Theorem~4, such a billiard system, observed with any finite-valued observation function, gives the same predictions as some observed semi-Markov process (disregarding differences in probabilistic predictions of less than $\varepsilon$).

Consequently, I conclude that the answer to the
question advanced at the beginning of this section is positive: \textit{there are deterministic systems in science which are observationally equivalent to stochastic processes in science at every observation level}. Compared to the results in Section~\ref{AI1} these results are \textit{stronger} because there is observational equivalence at every observation level (and not just for some observation).\footnote{The results here concern the same kind of descriptions as in Section~\ref{AI1}, namely continuous Bernoulli systems, semi-Markov processes and $n$-step semi-Markov processes.}
Note that the increase in strength is very different for the move from the first to the second results (observational equivalence not just to abstract mathematical descriptions but between descriptions in science) and for the move from the second to the third results (observational equivalence between descriptions in science at every, and not just at some, observation level). The third results imply that arguments that the evidence favours a deterministic description fail when they rely on the premise that fine-enough observations of deterministic systems in science yield probability distributions different from the ones of stochastic processes in science. And for this reason the argument advanced at the beginning of this subsection fails.
This, of course, is not to claim that there can be no sound arguments for preferring a deterministic or stochastic description, but just that arguments based on alleged general differences between probability distributions of deterministic and stochastic descriptions in science will not work.\footnote{
The question of which description is preferable is an interesting one but requires another paper. xxx is a paper devoted to this question.}

In this context note the following. Assume one wants to find out, by applying observation functions, whether
the evidence favours a deterministic or stochastic description. According to Theorem~3, certain deterministic systems are observationally equivalent at every observation level to $n$-step semi-Markov processes. Scientists use $n$-step semi-Markov processes in many contexts (from physics and biology to the social sciences) where they do not, and we think it is rational that they do not, automatically conclude that this indicates a deterministic description. Thus
it is plausible to argue that $n$-step semi-Markov processes do not indicate deterministic descriptions. Still, $n$-step semi-Markov processes show correlations in the following sense: the next outcome depends on the past $n$ outcomes, not only on the previous outcome. Suppose that, because of this, $n$-step semi-Markov processes are taken to indicate a deterministic description. Then one still cannot conclude that fine-enough observations of Bernoulli systems indicate deterministic descriptions. This is so because all one needs is that the observed probability distributions could derive from stochastic process showing no correlations. And, according to Theorem~4, all observations could derive from semi-Markov processes, where the next outcome  depends only on the past outcome and thus there are no such correlations.\footnote{It seems worth noting that for the semi-Markov processes of Theorem~4 any arbitrary outcome will be followed by at least two different outcomes. And for many Bernoulli systems any of these semi-Markov processes is such that none of the transition probabilities is close to $0$ or $1$ (see Subsection~\ref{P8}).}

All in all, the results on observational equivalence are stronger than what one might have expected: observational equivalence holds even between deterministic and stochastic descriptions \textit{found in science} at \textit{every observation level}.

\section{Previous Philosophical Discussion}\label{PM2}
There are hardly any philosophical reflections in the mathematics literature on observational equivalence. Because I discussed the continuous-time results, I can comment on the main exception:
\begin{quote}
Our theorem [Theorem~\ref{T5}] also tells us that certain semi-Markov systems could be thought of as being produced by Newton's laws (billiards seen through a deterministic viewer) or by coin-flipping. This may mean that there is no philosophical distinction between processes governed by roulette wheels and processes governed by Newton's laws. $\{$The popular literature emphasises the distinction between ``deterministic chaos'' and ``real randomness''.$\}$ In this connection we should note that our model for a stationary process (\text{\S} $1.2)$ [the deterministic representation] means that random processes have a deterministic model. This model, however, is abstract, and there is no
reason to believe that it can be endowed with any special additional
structure. Our point is that we are comparing, in a strong sense,
Newton's laws and coin flipping.\footnote{The text in braces is in a footnote.} (Ornstein and Weiss~1991,~39--40)
\end{quote}

It is hard to tell what this comment expresses because it is vague and unclear.\footnote{Ornstein and Weiss are mathematicians and not philosophers. So one should not blame them for misguided philosophical claims. Still, one needs to know whether their philosophical claims are tenable; thus I criticise them.}
For instance, why do Ornstein and Weiss highlight coin flipping even though Theorem~\ref{T5} does not tell us anything about coin flipping (Bernoulli processes) but only about semi-Markov processes?
Disregarding that, possibly, Ornstein and Weiss think that semi-Markov processes are random and hence claim that both deterministic and stochastic descriptions can be random. This is widely acknowledged in  philosophy (e.g., Eagle~2005).
Or maybe Ornstein and Weiss~(1991) want to say that deterministic systems in science, when observed with specific observation functions, can be observationally equivalent to stochastic processes in science or, if semi-Markov processes are random, even random processes.\footnote{If $(M,\Sigma_{M},\mu,T_{t})$ and a semi-Markov process $\{Z_{t};\,t\in\field{R}\}$ are $\varepsilon$-congruent, there is a finite-valued $\Phi$ such that $\{\Phi(T_{t});\,\, t\in\field{R}\}$ is $\{Z_{t};\,t\in\field{R}\}$ (cf.~Subsection~\ref{muede}).} This is true and important. However, as discussed in Section~\ref{AI1}, this was generally known before Theorem~\ref{T5} was proven and has been established by theorems which are weaker than Theorem~\ref{T5}. One might have expected Ornstein and Weiss~(1991) to say that Theorem~\ref{T5} shows what I argued that it does, namely that deterministic systems in science are observationally equivalent at every observation level to stochastic processes in science (cf.~Subsection~\ref{bermading}). But they do not seem to say this: because, if they did, it would be unclear why the deterministic representation is mentioned; also, they do not talk about all possible observation levels.

In any case, even if Theorem~\ref{T5} establishes observational equivalence, it is not true that ``this may mean that there is no philosophical distinction between processes governed by roulette wheels and processes governed by Newton's laws'' in the sense that there is no conceptual distinction between deterministic and indeterministic descriptions. Regardless of any results on observational equivalence, this distinction remains.

In the philosophy literature the significance of Theorem~\ref{T5} is taken to be that (at every observation level) there is a choice between a deterministic description in science and a stochastic process: ``What is fundamental is that [...]\ we are in a position to choose either between a deterministic or stochastic model'' (Suppes and de~Barros 1996,~196). ``The fact that a Bernoulli flow can be partitioned in such a way as to yield a (semi-) Markov process illustrates what has been acknowledged all along: Some deterministic systems, when partitioned, generate stochastic processes'' (Winnie 1998,~317).

Theorem~\ref{T5} indeed implies that there is a choice (at every observation level) between  deterministic descriptions in science and stochastic descriptions. Yet note that this is already shown by the first level (cf.~Subsection~\ref{DRS}).
For Theorem~\ref{epC} already shows that for every finite-valued observation function there is a choice between a nontrivial stochastic or a deterministic description. This implies that, according to our first notion, the deterministic systems to which Theorem~\ref{epC} applies are observationally equivalent at every observation level to nontrivial stochastic processes. Formally:
\begin{proposition}\label{Pfinal2}
Let $(M,\Sigma_{M},\mu,T_{t})$ be a measure-preserving deterministic system where $(M,d_{M})$ is separable and where $\Sigma_{M}$ contains all open sets of $(M,d_{M})$. Assume that $(M,\Sigma_{M},\mu,T_{t})$ satisfies the assumption of Theorem~\ref{epC} and has finite KS-entropy.\footnote{Deterministic systems in science have finite KS-entropy (Ornstein and Weiss~1991, 19).} Then for every $\varepsilon>0$ there is a  stochastic process $\{Z_{t};\,t\in\field{R}\}$ with outcome space $M_{O}=\cup_{l=1}^{h}o_{l}$, $h\in\field{N}$, such that $\{Z_{t};\,t\in\field{R}\}$ is $\varepsilon$-congruent to $(M,\Sigma_{M},\mu,T_{t})$, and for all $k\in\field{R}^{+}$ there are $o_{i},o_{j}\in M_{O}$ such that $0<P\{Z_{t+k}\!=\!o_{j}\,|\,Z_{t}\!=\!o_{i}\}<1$.\end{proposition}
This proposition is easy to establish (see Subsection~\ref{Pfinal2P}).
And, clearly, according to the second and third notion, every  deterministic system $(M,\Sigma_{M},\mu,T_{t})$ to which Theorem~\ref{epC} applies is observationally equivalent at every observation level to nontrivial stochastic processes. This is so because the second and third notion quantifies over all finite-valued observations $\Phi$, and $\{\Phi(T_{t}); t\in\field{R}\}$ is nontrivial by Theorem~\ref{epC}. Thus $\{\Phi(T_{t}); t\in\field{R}\}$ strongly and weakly $(\Phi,\varepsilon)$-simulates $(M,\Sigma_{M},\mu,T_{t})$ for every $\varepsilon>0$.

The significance of a theorem is constituted by the new knowledge added and by the conclusions shown by the theorem that have not already been shown by weaker theorems. As just explained, that there is a choice between deterministic and stochastic descriptions follows from much \textit{weaker} theorems than Theorem~\ref{T5} and was known to the community long before Theorem~\ref{T5} was proven. Consequently, \textit{the choice (at every observation level) between a deterministic description in science and a stochastic description cannot be the significance of Theorem~\ref{T5}.} Instead, its significance is that deterministic systems in science are observationally equivalent at every observation level even to stochastic processes of the kinds found in science.\footnote{Suppes~(1999,~182) and  Winnie (1998,~317) similarly claim that the philosophical significance of the result that some \textit{discrete-time} systems are $\varepsilon$-congruent for all $\varepsilon>0$ to \textit{discrete-time} Markov processes is that there is a choice between deterministic and stochastic descriptions. Werndl (2009a,~Section~4.2.2) already criticised that these claims are too weak but did not explain why (so I did this here). }

Moreover, Suppes and de~Barros (1996,~198--200) and Suppes~(1999,~189, 192) seem to think, wrongly, that what it means for a deterministic system to be $\varepsilon$-congruent to a certain type of stochastic process for every $\varepsilon>0$ (the first notion of observational equivalence) is that the system observed with any finite-valued observation function yields a process of a certain type: ``we can form a finite partition [...]\ of the phase space of possible trajectories'' (Suppes~1999,~189). That is, they wrongly think that the first notion of observational equivalence expresses something like my second notion.
 As discussed in Subsection~\ref{muede}, the first and the second notion are quite different (e.g.,~only the latter tells us what happens if any arbitrary finite-valued observation function is applied).

\section{Conclusion}
This paper presented and philosophically assessed three types of results on the observational equivalence of deterministic and indeterministic descriptions. These results were about continuous-time measure-theoretic deterministic systems and stochastic processes, which are ubiquitous in science. The main contribution of this paper was to show that there are results on observational equivalence for continuous-time  descriptions of increasing strength (namely, about descriptions in science and for every observation level) and to assess their philosophical significance and consequences.

The first results were about observational equivalence to abstract mathematical descriptions. I showed how, when starting with deterministic systems, one finds observationally equivalent stochastic processes, and conversely. The second results were stronger in the sense that they were not about observational equivalence to abstract mathematical descriptions but showed observational equivalence between descriptions of the types found in science. They establish, e.g., that billiard systems with convex obstacles are observationally equivalent to semi-Markov processes. Here I also discussed Kolmogorov's failed attempt to separate deterministic systems in science from deterministic systems which yield stochastic processes in science. Compared to the second results, the third results were stronger because they concerned observational equivalence between descriptions in science at every (and not just at some) observation level. I introduced two new meanings of `observational equivalence at every observation level'. And I showed that deterministic systems in science can be observationally equivalent even at every observation level to stochastic processes in science. For example, billiard systems with convex obstacles are observationally equivalent at every observation level to semi-Markov or $n$-step semi-Markov processes.

A consequence of these results is as follows. Suppose one wants to find out whether a phenomenon is best described as deterministic or stochastic. Then one might think of arguing that there is evidence for a deterministic or stochastic description by appealing to the difference between the probability distributions of stochastic processes in science and the probability distributions of (possibly fine-enough) observations of deterministic systems in science. The results of this paper showed that such arguments fail because there is no general difference between the respective probability distributions. Finally, I criticised the previous philosophical discussion of mathematicians and philosophers on observational equivalence, arguing that some of the claims of Ornstein and Weiss (1991), Suppes and de~Barros (1996), Suppes~(1999) and Winnie~(1998) are misleading.

\section{Appendix}
\subsection{\textit{Semi-Markov and $n$-step Semi-Markov Processes}}\label{A1}
First, we need the following definition.
\begin{definition}\label{DSP}
A discrete stochastic process $\{Z_{t};\,t\in\field{Z}\}$ is a family of random variables $Z_{t},\,\,t\in\field{Z}$, from a probability space $(\Omega,\Sigma_{\Omega},\nu)$ to a measurable space $(\bar{M},\Sigma_{\bar{M}})$.
\end{definition}
Semi-Markov and $n$-step semi-Markov processes are defined via Markov and $n$-step Markov processes (for a definition of these well-known processes, see Doob~(1953), where it is also explained what it means for these processes to be irreducible and aperiodic).

A semi-Markov process\footnote{The term `semi-Markov process' is not defined unambiguously. I follow Ornstein and Weiss~(1991).} is defined with help of a discrete stochastic process \linebreak[4] $\{(S_{k},T_{k}),\,k\in\field{Z}\}$.  $\{S_{k};\,\,k\in\field{Z}\}$ describes the successive outcomes $s_{i}$ visited by the semi-Markov process, where at time $0$ the outcome is $S_{0}$. $T_{0}$ is the time interval after which there is the first jump of the semi-Markov process after the time $0$, $T_{-1}$ is the time interval after which there is the last jump of the process before time $0$, and all other $T_{k}$ similarly describe the time-intervals between jumps of the process. Because at time $0$ the semi-Markov process is in $S_{0}$ and the process is in $S_{0}$ for the time $u(S_{0})$, $T_{-1}=u(S_{0})-T_{0}$.

Technically, $\{(S_{k},T_{k}),\,k\in\field{Z}\}$ satisfies the following conditions: (i) $S_{k}\!\in\! S\!=\!\{s_{1},\ldots,s_{N}\}$, $N\in\field{N}$; $T_{k}\in U=\{u_{1},\ldots,u_{\bar{N}}\}$, $\bar{N}\in\field{N}$, $\bar{N}\leq N$, for $k\neq 0,-1$, where $u_{i}\in\field{R}^{+}$, $1\leq i \leq \bar{N}$; $T_{0}\in (0,u(S_{0})]$, $T_{-1}\in [0,u(S_{0}))$, where $u:S\rightarrow U, s_{i}\rightarrow u(s_{i})$, is surjective; and hence $\bar{M}=S\times [0, \max_{i}u_{i}]$; (ii) $\Sigma_{\bar{M}}=\field{P}(S)\times L([0,\max_{i}u_{i}])$, where $L([0,\max_{i}u_{i}])$ is the Lebesgue $\sigma$-algebra on $[0,\max_{i}u_{i}]$; (iii) $\{S_{k};\,\,k\in\field{Z}\}$ is a stationary irreducible and aperiodic Markov process with outcome space $S$; $p_{s_{i}}=P\{S_{0}=s_{i}\}>0$, for all $i$, $1\leq i\leq N$; (iv) $T_{k}=u(S_{k})$ for $k\geq 1$, $T_{k}=u(S_{k-1})$ for $k\leq -2$, and $T_{-1}=u(S_{0})-T_{0}$; (v) for all $i$, $1\leq i \leq N$,
$P(T_{0}\in A\,|\,S_{0}=s_{i})=\int_{A} 1/u(s_{i})d\lambda$ for all $A\in L((0,u(s_{i})])$, where $L((0,u(s_{i})])$ is the Lebesgue $\sigma$-algebra and $\lambda$ is the Lebesgue measure on $(0,u(s_{i})]$.

\begin{definition}
A $\mathrm{semi}$-$\mathrm{Markov}$ $\mathrm{process}$ is a process  $\{Z_{t};\,t\in\field{R}\}$ with outcome space $S$ constructed via a process $\{(S_{k},T_{k}),\,k\in\field{Z}\}$ as follows:
\begin{eqnarray}
Z_{t}&\!\!\!\!=\!\!\!\!&S_{0}\,\,\textnormal{for}\,\,-T_{-1}\leq t<T_{0},\nonumber\\
Z_{t}&\!\!\!\!=\!\!\!\!&S_{k}\,\,\textnormal{for}\,\,T_{0}+\ldots + T_{k-1}\leq t < T_{0}+\ldots + T_{k}; k\geq 1\,\,\textnormal{and thus}\,\,t\geq T_{0},\nonumber \\
Z_{t}&\!\!\!\!=\!\!\!\!&S_{-k}\,\textnormal{for}\!-\!T_{-1}\!-\!\ldots\!-\!T_{-k-1}\leq t < \!-\!T_{-1}\!-\!\ldots \!-\!T_{-k}; k\geq 1\,\textnormal{and thus}\,t<\!-\!T_{-1},\nonumber
\end{eqnarray}
and for all $i$, $1\leq i \leq N$,
\begin{equation}
P(Z_{0}=s_{i})=\frac{p_{s_{i}}u(s_{i})}{p_{s_{1}}u(s_{1})+\ldots +p_{s_{N}}u(s_{N})}.
\end{equation}
\end{definition}

Semi-Markov processes are stationary (Ornstein~1970 and~1974,~56--61). In this paper I assume that the elements of $U$ are irrationally related ($u_{i}$ and $u_{j}$ are \textit{irrationally related} iff $\frac{u_{i}}{u_{j}}$ is irrational; and the elements of $U=\{u_{1},\ldots,u_{\bar{N}}\}$ are \textit{irrationally related} iff for all $i,j$, $i\neq j$, $u_{i}$ and $u_{j}$ are irrationally related).

$n$-$step$ semi-Markov processes\footnote{Again, the term `$n$-step semi-Markov process' is not used unambiguously; I follow Ornstein and Weiss (1991).} generalise semi-Markov processes.
\begin{definition}
$n$-$\mathrm{step}$ $\mathrm{semi}$-$\mathrm{Markov}$ $\mathrm{processes}$ are defined like semi-Markov processes except that  condition (iii) is replaced by: (iii') $\{S_{k};\,\,k\in\field{Z}\}$ is a stationary irreducible and aperiodic $n$-step Markov process with outcome space $S$ and $p_{s_{i}}=P\{S_{0}=s_{i}\}>0$, for all $i$, $1\leq i\leq N$.
\end{definition}
Again, $n$-step semi-Markov processes are stationary (Park~1982), and I assume that the elements of $U$ are irrationally related.

\subsection{\textit{Proof of Theorem~\ref{epC}}}\label{A2}
\textit{Iff for a measure-preserving deterministic system $(M,\Sigma_{M},\mu,T_{t})$ there does not exist a $n\in\field{R}^{+}$ and a $C\in\Sigma_{M}$, $0<\mu(C)<1,$ such that, except for a set of measure zero (esmz.) $T_{n}(C)=C$, then the following holds: for every nontrivial finite-valued observation function $\Phi:M\rightarrow M_{O}$, every $k\in\field{R}^{+}$ and $\{Z_{t};\,t\in\field{R}\}=\{\Phi(T_{t});\,t\in\field{R}\}$ there are $o_{i},o_{j}\in M_{O}$ with $0<P\{Z_{t+k}\!=\!o_{j}\,|\,Z_{t}\!=\!o_{i}\}<1$.}\\

\noindent \textit{Proof}:
We need the following definitions.
\begin{definition}\label{DS}
A $\mathrm{discrete}$ $\mathrm{measure}$-$\mathrm{preserving}$ $\mathrm{deterministic}$ $\mathrm{system}$ is a quadruple \linebreak[4] $(M,\Sigma_{M},\mu,T)$ where $(M,\Sigma_{M},\mu)$ is a probability space and $T:M\rightarrow M$ is a bijective measurable function such that $T^{-1}$ is measurable and $\mu(T(A))=\mu(A)$ for all $A\in\Sigma_{M}$.
\end{definition}
\begin{definition} \label{ergodic} A discrete measure-preserving deterministic system $(M,\Sigma_{M},\mu,T)$ is $\mathrm{ergodic}$ iff there is no $A\in\Sigma_{M}$, $0<\mu(A)<1$, such that, esmz., $T(A)=A$.
\end{definition}
$(M,\Sigma_{M},\mu,T)$ is ergodic iff for all $A,B\in\Sigma_{M}$ (Cornfeld et al.~1982,~14--15): \begin{equation}\label{erg}
\lim_{n\rightarrow \infty}\frac{1}{n}\sum_{i=1}^{n}(\mu(T^{n}(A)\cap B)-\mu(A)\mu(B))=0.
\end{equation}
\begin{definition}\label{partition}
$\alpha=\{\alpha_{1},\ldots,\alpha_{n}\}$, $n\in\field{N}$, is a $\mathrm{partition}$ of $(M,\Sigma_{M},\mu)$ iff $\alpha_{i}\in\Sigma_{M},\,\,\mu(\alpha_{i})>0,$ for all $i$, $1\leq i\leq n$, $\alpha_{i}\cap\alpha_{j}=\emptyset$ for all $i\neq j$, $1\leq i,j\leq n$, and $M=\bigcup_{i=1}^{n}\alpha_{i}$.\end{definition} A partition is nontrivial iff $n\geq 2$.
Given two partitions $\alpha=\{\alpha_{1},\ldots,\alpha_{n}\}$ and $\beta=\{\beta_{1},\ldots,\beta_{l}\}$ of $(M,\Sigma_{M},\mu)$, $\alpha\vee\beta$ is the partition $\{\alpha_{i}\cap\beta_{j}\,|\,i=1,\ldots,n;j=1,\ldots,l\}$. Given a deterministic system $(M,\Sigma_{M},\mu,T_{t})$, if $\alpha$ is a partition, $T_{t}\alpha=\{T_{t}(\alpha_{1}),\ldots,T_{t}(\alpha_{n})\}$, $t\in\field{R}$, is also a partition.

Note that any finite-valued observation function $\Phi$ (cf.~Subsection~\ref{Det}) can be written as: $\Phi(m)=\sum_{i=1}^{n}o_{i}\chi_{\alpha_{i}}(m)$,
$M_{O}=\{o_{i}\,|\,1\leq i\leq n\}$, for some partition $\alpha$
of $(M,\Sigma_{M},\mu)$, where $\chi_{A}$
is the characteristic function of $A$ (cf.~Cornfeld et al.~1982,~179).\footnote{That is, $\chi_{A}(m)=1$ for $m\in A$ and $0$ otherwise.} It suffices to prove the following:
\begin{quote}$(*)$ Iff for $(M,\Sigma_{M},\mu,T_{t})$ there does not exist an $n\in\field{R}^{+}$ and a
$C\in\Sigma_{M}$, $0<\mu(C)<1,$ such that, esmz., $T_{n}(C)=C$, then the following holds: for any nontrivial partition $\alpha=\{\alpha_{1},\ldots,\alpha_{r}\}$, $r\in\field{N}$, and all $k\in\field{R}^{+}$ there is an $i\in\{1,\ldots,r\}$ such that for all $j,\,1\!\leq\! j\!\leq\! r$,  $\mu(T_{k}(\alpha_{i})\!\setminus\!\alpha_{j})\!>\!0$.
\end{quote}
For finite-valued observation functions  are of the form $\sum_{l=1}^{r}o_{l}\chi_{\alpha_{l}}(m)$, where $\alpha=\{\alpha_{1},\ldots,\alpha_{r}\}$ is a partition and $M_{O}=\cup_{l=1}^{r}o_{l}$. Consequently, the right hand side of $(*)$ expresses that for any nontrivial finite-valued $\Phi:M\rightarrow M_{O}$ and all $k\in\field{R}^{+}$ there is an $o_{i}\in M_{O}$  such that for all $o_{j}\in M_{O}$,  $P\{Z_{t+k}=o_{j}\,|\,Z_{t}=o_{i}\}<1$, or equivalently, that there are $o_{i},\,o_{j}\in M_{O}$ with $0<P\{Z_{t+k}\!=\!o_{j}\,|\,Z_{t}\!=\!o_{i}\}<1$.\\

\noindent $\Leftarrow:$ Assume that there is an $n\in\field{R}^{+}$ and a $C\in\Sigma_{M}$, $0<\mu(C)<1$, such that, esmz., $T_{n}(C)=C$. Then for  $\alpha=\{C,M\setminus C\}$ it holds that $\mu(T_{n}(C)\setminus C)=0$ and $\mu(T_{n}(M\setminus C)\setminus (M\setminus C))=0$.\\

\noindent $\Rightarrow:$ Assume that the conclusion of $(*)$ does not hold, and hence that there is a nontrivial partition $\alpha$ and a $k\in\field{R}^{+}$ such that for each $\alpha_{i}$ there is an $\alpha_{j}$ with, esmz.,
$T_{k}(\alpha_{i})\subseteq \alpha_{j}$. From the assumptions it follows that for every $k\in\field{R}^{+}$ the discrete deterministic system $(M,\Sigma_{M},\mu,T_{k})$ is ergodic (cf.~Definition~\ref{ergodic}).

\textit{Case 1}: For every $i$ there is a $j$ such that, esmz., $T_{k}(\alpha_{i})= \alpha_{j}$. Because the discrete system $(M,\Sigma_{M},\mu,T_{k})$ is ergodic (equation (\ref{erg})), there is an $h\in\field{N}$ such that, esmz., $T_{kh}(\alpha_{1})=\alpha_{1}$. But this is in contradiction with the assumption that it is not the case that there exists an $n\in\field{R}^{+}$ and a
$C\in\Sigma_{M},$ $0<\mu(C)<1,$ such that, esmz., $T_{n}(C)=C$.

\textit{Case 2}: There exists an $i$ and a $j$ with, esmz.,
$T_{k}(\alpha_{i})\subset \alpha_{j}$ and with $\mu(\alpha_{i})<\mu(\alpha_{j})$. Because the discrete system $(M,\Sigma_{M},\mu,T_{k})$ is ergodic (equation (\ref{erg})), there is a $h\in\field{N}$ such that, esmz., $T_{hk}(\alpha_{j})\subseteq\alpha_{i}$. Hence $\mu(\alpha_{j})\leq\mu(\alpha_{i})$, contradicting $\mu(\alpha_{i})<\mu(\alpha_{j})\leq\mu(\alpha_{i})$.

\subsection{\textit{Proof of Theorem~\ref{semiMarkov1}}}\label{P8}
\textbf{Theorem~\ref{semiMarkov1}} \textit{Let $(M,\Sigma_{M},\mu,T_{t})$ be a continuous Bernoulli system. Then for every finite-valued $\Phi$ and every $\varepsilon>0$ a semi-Markov process $\{Z_{t},\,t\in\field{R}\}$ weakly  $(\Phi,\varepsilon)$-simulates $(M,\Sigma_{M},\mu,T_{t})$.}\\

\noindent \textit{Proof}: Let $(M,\Sigma_{M},\mu,T_{t})$ be a continuous Bernoulli system, let $\Phi:M\rightarrow S, S=\{s_{1},\ldots,s_{N}\}$, $N\in\field{N}$, be a surjective observation function, and let $\varepsilon>0$. Theorem~\ref{T6} implies that there is a surjective observation function $\Theta:M\rightarrow S,\,\,\Theta(m)=\sum_{i=1}^{N}s_{i}\chi_{\alpha_{i}}(m)$, for a partition $\alpha$ (cf.~Definition~\ref{partition}), such that $\{Y_{t}=\Theta(T_{t});\,t\in\field{R}\}$ is an $n$-step semi-Markov process with outcomes $s_{i}$ and corresponding times $u(s_{i})$ which strongly $(\Phi,\varepsilon)$-simulates $(M,\Sigma_{M},\mu,T_{t})$.

I need the following definition:
\begin{definition}\label{factor}
$(M_{2},\Sigma_{M_{2}},\mu_{2},T^{2}_{t})$ is a $\mathrm{factor}$ of  $(M_{1},\Sigma_{M_{1}},\mu_{1},T^{1}_{t})$ (where both deterministic systems measure-preserving) iff there are
 $\hat{M}_{i}\subseteq M_{i}$ with
$\mu_{i}(M_{i}\setminus\hat{M}_{i})=0$,
$T^{i}_{t}\hat{M}_{i}\subseteq\hat{M}_{i}$ for all $t\,\,(i=1,2)$, and there is a function $\phi:\hat{M}_{1}\!\rightarrow \!\hat{M}_{2}$ such that (i) $\phi^{-1}(B)\!\in\!\Sigma_{M_{1}}$ for all
$B\!\in\!\Sigma_{M_{2}},A\subseteq\hat{M}_{2}$;
(ii) $\mu_{1}(\phi^{-1}(B))=\mu_{2}(B)$ for all
$B\in\Sigma_{M_{2}},\,B\subseteq\hat{M}_{2}$; (iii)
$\phi(T^{1}_{t}(m))=T^{2}_{t}(\phi(m))$ for all $m\in\hat{M}_{1},\,t\in\field{R}$.
\end{definition}
Note that the deterministic representation $(X,\Sigma_{X},\mu_{X},W_{t},\Lambda_{0})$ of $\{Y_{t};\,\,t\in\field{R}\}$ is a factor of $(M,\Sigma_{M},\mu,T_{t})$ via $\phi(m)=r_{m}$ ($r_{m}$ is the realisation of $m$) (cf.~Ornstein and Weiss,~1991,~18).

Now I construct a measure-preserving system $(K,\Sigma_{K},\mu_{K},R_{t})$ as follows. Let \linebreak[4] $(\Omega,\Sigma_{\Omega},\mu_{\Omega},V,\Xi_{0})$, $\Xi_{0}(\omega)=\sum_{i=1}^{N}s_{i}\chi_{\beta_{i}}(\omega)$, where $\beta$ is a partition, be the deterministic representation of $\{S_{k};\,k\in\field{Z}\}$, the irreducible and aperiodic $n$-step Markov process corresponding to $\{Y_{t};\,t\in\field{R}\}$.
Let $f:\Omega\rightarrow\{u_{1},\ldots,u_{N}\}$, $f(\omega)=u(\Xi_{0}(\omega))$.
Define $K$ as $\cup_{i=1}^{N}K_{i}=\cup_{i=1}^{N}(\beta_{i}\times [0,u(s_{i})))$. Let $\Sigma_{K_{i}}$, $1\leq i \leq N$, be the product $\sigma$-algebra $(\Sigma_{\Omega}\cap\beta_{i})\times L([0,u(s_{i})))$ where $L([0,u(s_{i})))$ is the Lebesgue $\sigma$-algebra of $[0,u(s_{i}))$. Let $\mu_{K_{i}}$ be the product measure
\begin{equation} (\mu_{\Omega}^{\Sigma_{\Omega}\cap\beta_{i}}\times\lambda([0,u(s_{i}))))/\sum_{j=1}^{N}u(s_{j})\mu_{\Omega}(\beta_{j}),
\end{equation}
 where $\lambda([0,u(s_{i})))$ is the Lebesgue measure on $[0,u(s_{i}))$ and $\mu_{\Omega}^{\Sigma_{\Omega}\cap\beta_{i}}$ is the measure $\mu_{\Omega}$ restricted to $\Sigma_{\Omega}\cap\beta_{i}$. Let $\Sigma_{K}$ be the $\sigma$-algebra generated by $\cup_{i=1}^{N}\Sigma_{K_{i}}$. Define a pre-measure $\bar{\mu}_{K}$ on $H=(\cup_{i=1}^{N}(\Sigma_{\Omega}\cap\beta_{i}\times L([0,s_{i}))))\cup K$ by $\bar{\mu}_{K}(K)=1$ and $\bar{\mu}_{K}(A)=\mu_{K_{i}}(A)$ for $A\in\Sigma_{K_{i}}$, and let $\mu_{K}$ be the unique extension of this pre-measure to a measure on $\Sigma_{K}$.
Finally, $R_{t}$ is defined as follows: let the state of the system at time zero be  $(k,v)\in K$, $k\in\Omega,v<f(k)$; the state moves vertically with unit velocity, and just before it reaches $(k,f(k))$ it jumps to $(V(k),0)$ at time $f(k)-v$; then it again moves vertically with unit velocity,
and just before it reaches $(V(k),f(V(k))))$ it jumps to $(V^{2}(k),0)$ at time $f(V(k))+f(k)-v$, and so on. $(K,\Sigma_{K},\mu_{K},R_{t})$ is a measure-preserving system (a `flow built under the function $f$'). $(X,\Sigma_{X},\mu_{X},W_{t})$ is proven to be isomorphic (via a function $\psi$) to $(K,\Sigma_{K},\mu_{K},R_{t})$ (Park~1982).

Consider $\gamma=\{\gamma_{1},\ldots,\gamma_{l}\}\!=\!\beta\vee V\beta\vee\ldots\vee V^{n-1}\beta$ and $\Pi(\omega)=\sum_{j=1}^{l}o_{j}\chi_{\gamma_{j}}(\omega)$, $o_{i}\neq o_{j}$ for $i\neq j,\,1\leq i,\,j\leq l$.
I now show that $\{B_{t}=\Pi(V^{t}(\omega));\,\,t\in\field{Z}\}$ is an irreducible and aperiodic Markov process.
 By construction, for all $t$ and all $i$, $1\leq i\leq l$, there are $q_{i,0},\ldots, q_{i,n-1}\in S$
such that
\begin{equation}
P\{B_{t}=o_{i}\}=P\{S_{t}=q_{i,0}, S_{t+1}=q_{i,1},\ldots,S_{t+n-1}=q_{i,n-1}\}.
\end{equation}
Therefore, for all $k\in\field{N}$ and all $i,j_{1},\ldots,j_{k},\,\,1\leq i,j_{1},\ldots,j_{k}\leq l$:
\begin{equation}
P\{B_{t+1}=o_{i}\,|\,B_{t}=o_{j_{1}},\ldots,B_{t-k+1}=o_{j_{k}}\}=
\end{equation}
\begin{eqnarray}
P\{S_{t+1}\!\!=\!\!q_{i,0},\!...,S_{t+n}\!\!=\!\!q_{i,n-1}|S_{t}\!\!=\!\!q_{j_{1},0},\!...,S_{t+n-1}\!\!=\!\!q_{i,n-2},\!S_{t-1}\!\!=\!\!q_{j_{2},0},\!...\!,S_{t-k+1}\!\!=\!\!q_{j_{k},0}\}\nonumber \\
=P\{S_{t+1}=q_{i,0},\ldots,S_{t+n}=q_{i,n-1}\,|\,S_{t}=q_{j_{1},0},\ldots,S_{t+n-1}=q_{i,n-2}\}\nonumber\\
=P\{B_{t+1}=o_{i}\,|\,B_{t}=o_{j_{1}}\},\nonumber
\end{eqnarray}
if $P\{B_{t+1}=o_{i},B_{t}=o_{j_{1}},\ldots,B_{t-k+1}=o_{j_{k}}\}>0$.
Hence $\{B_{t};\,t\in\field{Z}\}$ is a Markov process. Now a discrete measure-preserving deterministic system  $(M,\Sigma_{M},\mu,T)$ is mixing iff for all $A,B\in\Sigma_{M}$
\begin{equation}
\lim_{t\rightarrow\infty}\mu(T^{t}(A)\cap B)=\mu(A)\mu(B).
\end{equation}
The deterministic representation of every irreducible and aperiodic $n$-step Markov process, and hence $(\Omega,\Sigma_{\Omega},\mu_{\Omega},V,\Xi_{0})$, is mixing (Ornstein~1974,~45--47). This implies that $\{B_{t};\,t\in\field{Z}\}$ is irreducible and aperiodic.

Let $\Delta(k)=\sum_{i=1}^{l}o_{i}\chi_{\gamma_{i}\times[0,u(o_{i}))}(k)$,  where $u(o_{i}),\,1\leq i\leq l$, is defined as follows: $u(o_{i})=u(s_{r})$ where $\gamma_{i}\subseteq \beta_{r}$. It follows immediately that $\{X_{t}=\Delta(R_{t});\,t\in\field{R}\}$ is a semi-Markov process. Consider the surjective function $\Psi:M\rightarrow\{o_{1},\ldots,o_{l}\}$, $\Psi(m)=\Delta(\psi(\phi(m)))$ for $m\in\hat{M}$ and $o_{1}$ otherwise. Recall that $(X,\Sigma_{X},\mu_{X},W_{t})$ is a factor (via $\phi$) of $(M,\Sigma_{M},\mu,T_{t})$ and that $(X,\Sigma_{X},\mu_{X},W_{t})$ is isomorphic (via $\psi$) to $(K,\Sigma_{K},\mu_{K},R_{t})$. Therefore, $\{Z_{t}=\Psi(T_{t});\,t\in\field{R}\}$ is a semi-Markov process with outcomes $o_{i}$ and times $u(o_{i}),\,\,1\leq i \leq l$.

Now consider the surjective observation function $\Gamma:\{o_{1},\ldots,o_{l}\}\rightarrow S$, where $\Gamma(o_{i})=s_{r}$ for $\gamma_{i}\subseteq\beta_{r}$, $1\leq i\leq l$.
By construction, esmz., $\Gamma(\Psi(T_{t}(m)))=\Theta(T_{t}(m))=Y_{t}(m)$ for all $t\in\field{R}$. Hence, because $\{Y_{t};\,t\in\field{R}\}$ strongly $(\Phi,\varepsilon)$-simulates $(M,\Sigma_{M},\mu,T_{t})$, $\mu(\{m\in M\,|\,\Gamma(\Psi(m))\neq \Phi(m)\})<\varepsilon$.

\subsection{\textit{Proof of Proposition~\ref{Pfinal2}}}\label{Pfinal2P}
\textbf{Proposition~\ref{Pfinal2}} \textit{
Let $(M,\Sigma_{M},\mu,T_{t})$ be a measure-preserving deterministic system where $(M,d_{M})$ is separable and where $\Sigma_{M}$ contains all open sets of $(M,d_{M})$. Assume that $(M,\Sigma_{M},\mu,T_{t})$ satisfies the assumption of Theorem~\ref{epC} and has finite KS-entropy. Then for every $\varepsilon>0$ there is a  stochastic process $\{Z_{t};\,t\in\field{R}\}$ with outcome space $M_{O}=\cup_{l=1}^{h}o_{l}$, $h\in\field{N}$, such that $\{Z_{t};\,t\in\field{R}\}$ is $\varepsilon$-congruent to $(M,\Sigma_{M},\mu,T_{t})$, and for all $k\in\field{R}^{+}$ there are $o_{i},o_{j}\in M_{O}$ such that $0<P\{Z_{t+k}\!=\!o_{j}\,|\,Z_{t}\!=\!o_{i}\}<1$.}\\

\noindent \textit{Proof}:  A partition  $\alpha$ (cf.~Definition~\ref{partition}) is \textit{generating} for the discrete measure-preserving system $(M,\Sigma_{M},\mu,T)$ iff
for every $A\in\Sigma_{M}$ there is an $n\in\field{N}$ and a set $C$ of unions of elements in $\vee_{j=-n}^{n}T^{j}(\alpha)$ such that $\mu((A\setminus C)\cup(C\setminus A))<\varepsilon$ (Petersen~1983,~244). $\alpha$ is \textit{generating} for the measure-preserving system $(M,\Sigma_{M},\mu,T_{t})$ iff for all $A\in\Sigma_{M}$ there is a $\tau\in\field{R}^{+}$ and a set $C$ of unions of elements in $\bigcup_{\textnormal{all}\,\,m}\bigcap_{t=-\tau}^{\tau}(T^{-t}(\alpha(T^{t}(m))))$ such that $\mu((A\setminus C)\cup(C\setminus A))<\varepsilon$ ($\alpha(m)$ is the set $\alpha_{j}\in\alpha$ with $m\in\alpha_{j}$).

By assumption, there is a $t_{0}\in\field{R}^{+}$ such that the discrete system $(M,\Sigma_{M},\mu,T_{t_{0}})$ is ergodic (cf.~Definition~\ref{ergodic}). For discrete ergodic systems Krieger's~(1970) theorem implies that there is a partition $\alpha$ which is generating for $(M,\Sigma_{M},\mu,T_{t_{0}})$ and hence generating for $(M,\Sigma_{M},\mu,T_{t})$. Since $(M,d_{M})$ is separable, for every $\varepsilon>0$ there is a $r\in\field{N}$ and $m_{i}\in M$, $1\leq i \leq r$, such that $\mu(M\setminus\cup_{i=1}^{r}B(m_{i},\frac{\varepsilon}{2}))<\frac{\varepsilon}{2}$.
Because $\alpha$ is generating for $(M,\Sigma_{M},\mu,T_{t_{0}})$, for each $B(m_{i},\frac{\varepsilon}{2})$ there is an $n_{i}\in\field{N}$ and a $C_{i}$ of union of elements in $\vee_{j=-n_{i}}^{n_{i}}T_{jt_{0}}(\alpha)$ such that $\mu((B(m_{i},\frac{\varepsilon}{2})\setminus C_{i})\cup
(C_{i}\setminus B(m_{i},\frac{\varepsilon}{2})))<\frac{\varepsilon}{2r}$. Let $n\!=\!\max\{n_{i}\},\,\,\beta\!=\!\{\beta_{1},\ldots,\beta_{l}\}\!=\!\vee_{j=-n}^{n}T_{jt_{0}}(\alpha)$ and  $\Psi\!(m)=\!\sum_{i=1}^{l}o_{i}\chi_{\beta_{i}}(m)$ with $o_{i}\in\beta_{i}$.
Since $\Psi$ is a finite-valued, Theorem~\ref{epC} implies that for the process $\{\Psi(T_{t});\,\,t\in\field{R}\}$ for all $k\in\field{R}^{+}$ there are $o_{i},\,o_{j}$ such that $0<P\{Z_{t+k}\!=\!o_{j}\,|\,Z_{t}\!=\!o_{i}\}<1$. Because $\alpha$ is generating for $(M,\Sigma_{M},\mu,T_{t})$, $\beta$ is generating too. This implies that $(M,\Sigma_{M},\mu,T_{t})$ is isomorphic (via a function $\phi$) to the deterministic representation $(M_{2},\Sigma_{M_{2}},\mu_{2},T^{2}_{t},\Phi_{0})$ of $\{Z_{t};\,\,t\in\field{R}\}$ (Petersen~1983,~274). And, by construction, $d_{M}(m,\Phi_{0}(\phi(m)))<\varepsilon$ except for a set in $M$ smaller than $\varepsilon$.


\newpage
\section*{References}
\addcontentsline{toc}{section}{References}
\begin{list}{}{%
    \setlength{\labelwidth}{0pt}
    \setlength{\labelsep}{0pt}
    \setlength{\leftmargin}{24pt}
    \setlength{\itemindent}{-24pt}
  }

\item Berkovitz, Joseph, Roman Frigg and Fred Kronz 2006, ``The Ergodic Hierarchy, Randomness and {Hamiltonian}
  Chaos'', Studies in History and Philosophy of Modern Physics   37:661--691.

\item Butterfield, Jeremy  2005, ``Determinism and Indeterminism'',
\newblock Routledge Encyclopaedia of Philosophy Online.

\item Cornfeld, Isaak P., Sergej V. Fomin and Yakuv G. Sinai 1982, Ergodic Theory, Berlin: Springer.

\item Doob, Joseph L.  1953, Stochastic
  Processes, New York: John Wiley \& Sons.

\item Eagle, Antony  2005, ``Randomness Is
  Unpredictability'', The British Journal for the Philosophy of Science 56:749--790.

\item Eckmann, Jean-Paul and David Ruelle  1985, ``Ergodic Theory of Chaos and Strange Attractors'', Reviews of Modern
  Physics 57:617--654.

\item Feldman, Joel and Meir Smorodinksy  1971, ``Bernoulli Flows With Infinite Entropy'', The Annals of Mathematical Statistics 42:381--382.

\item Halmos, Paul  1944, ``In General a
  Measure-preserving Transformation is Mixing'', The Annals of Mathematics 45:786--792.

\item Halmos, Paul  1949, ``Measurable
  transformations'', Bulletin of the American Mathematical Society
  55:1015--1043.

\item Hopf, Eberhard  1932, ``Proof of {Gibbs'}
  Hypothesis on the Tendency Toward Statistical Equilibrium'', Proceedings
  of the National Academy of Sciences of the United States of America   18:333--340.

\item Janssen, Jacques and Nikolaos Limnios  1999,
  Semi-Markov Models and Applications, Dotrecht: Kluwer Academic Publishers.

\item Krieger, Wolfgang  1970, ``On Entropy and Generators of Measure-preserving Transformations'', Transactions of the  American Mathematical Society 149:453--456.

\item Lorenz, Edward 1963, ``Deterministic Nonperiodic
  Flow'', Journal of the Atmospheric Sciences 20:~130--141.

\item Luzzatto, Stefano, Ian Melbourne and Frederic Paccaut  2005, ``The {Lorenz} Attractor Is Mixing'',
  Communications in Mathematical Physics 260:393--401.

\item Ornstein, Donald 1970, ``Imbedding {Bernoulli} Shifts in Flows'', in Albrecht Dold and Beno Eckmann (eds.), Contributions to Ergodic Theory and Probability, Proceedings of the First Midwestern Conference on Ergodic Theory. Springer: Berlin, 178--218.

\item Ornstein, Donald 1974, Ergodic Theory,
  Randomness, and Dynamical Systems, New Haven and London: Yale University Press.

\item Ornstein, Donald and Giovanni Galavotti  1974, ``Billiards and {Bernoulli} Schemes'', Communications in Mathematical
  Physics 38:83--101.

\item Ornstein, Donald and Benjamin Weiss  1991,
  ``Statistical Properties of Chaotic Systems'', Bulletin of the American
  Mathematical Society 24:11--116.

\item Park, Kong 1982, ``A Special Family of Ergodic Flows and Their $\bar{d}$-Limits'', Israel Journal of Mathematics 42:343--353.

\item Petersen, K.  1983, {\em Ergodic Theory}, Cambridge: Cambridge University Press.

\item Radunskaya, Amy  1992, Statistical Properties
  of Deterministic Bernoulli Flows, Ph.D. Dissertation, Stanford: University of Stanford.

\item Sim\'{a}nyi, N\'{a}ndor  2003, ``Proof of the Boltzmann-Sinai Ergodic Hypothesis for Typical Hard Disk Systems'', Inventiones Mathematicae 154:~123--178.

\item Sinai, Yakov G.  1989, ``Kolmogorov's Work on
  Ergodic Theory'', The Annals of Probability 17:~833--839.

\item Strogatz, Steven H.  1994, Nonlinear Dynamics
  and Chaos, with Applications to Physics, Biology, Chemistry, and
  Engineering, New York: Addison Wesley.

\item Suppes, Patrick  1999, ``The Noninvariance of
  Deterministic Causal Models'', Synthese 121:~181--198.

\item Suppes, Patrick and Acacio de~Barros 1996,
  ``Photons, Billiards and Chaos'', in Paul Weingartner and Gerhard Schurz (eds.),
  Law and Prediction in the Light of Chaos Research. Berlin: Springer, 189--207.

\item Uffink, Jos  2007, ``Compendium to the
  Foundations of Classical Statistical Physics'', in Jeremy Butterfield
  and John Earman (eds.), Philosophy of Physics (Handbooks of the
  Philosophy of Science B). Amsterdam: North-Holland, 923--1074.

\item Werndl, Charlotte  2009a, ``Are Deterministic  Descriptions and Indeterministic Descriptions Observationally Equivalent?'', Studies in History and Philosophy of Modern Physics 40:232--242.

\item Werndl, Charlotte  2009b, ``What Are the New Implications of Chaos for Unpredictability?'', The British Journal for the Philosophy of Science 60:195--220.


\item Winnie, John 1998, ``Deterministic Chaos and
  The Nature of Chance'', in John Earman and John Norton (eds.), The
  Cosmos of Science -- Essays of Exploration. Pittsburgh: Pittsburgh University Press, 299--324.

\end{list}
\end{document}